\newif\ifJOC

\JOCfalse 

\ifJOC
\documentclass[ijoc,nonblindrev]{informs3}

\OneAndAHalfSpacedXII 

\usepackage{natbib}
 \bibpunct[, ]{(}{)}{,}{a}{}{,}%

\TheoremsNumberedThrough
\EquationsNumberedThrough  

\else 
\documentclass[final,3p]{scrartcl}
\usepackage{amstext}
\usepackage{amsmath,amssymb,amsthm,graphicx}

\usepackage[numbers]{natbib}
\setcitestyle{numbers}

\DeclareMathOperator*{\argmax}{arg\,max}
\DeclareMathOperator*{\argmin}{arg\,min}

\newtheorem{theorem}{Theorem}[section]
\newtheorem{definition}[theorem]{Definition}
\newtheorem{lemma}[theorem]{Lemma}
\newtheorem{corollary}[theorem]{Corollary}

\newtheorem{proposition}[theorem]{Proposition}

\theoremstyle{remark}
\newtheorem{remarkx}[theorem]{Remark}
\newenvironment{remark}
    {\pushQED{\qed}\remarkx}
    {\popQED\endremarkx}
    
\fi 

\usepackage{changepage}
\usepackage[utf8x]{inputenc}
\usepackage[T1]{fontenc}
\usepackage{lmodern}
\usepackage[english]{babel}
\usepackage{hyperref}
\usepackage{todonotes}
\usepackage{algpseudocode}
\usepackage[ruled,vlined]{algorithm2e}
\usepackage{subcaption}
\usepackage{pgfplots}
\usepackage{booktabs}
\usepackage{tabularx}
\usepackage{ltablex}
\usepackage{longtable}
\usepackage{comment}
\usepackage{url}
\usepackage{bbm}
\usepackage[normalem]{ulem}
\pgfplotsset{compat=1.9}

\newcommand{\bftab}{\fontseries{b}\selectfont} 
\newcommand{\inprod}[2]{{\langle #1,#2 \rangle}} 

\newcommand{\Acal}{\mathcal{A}}
\newcommand{\Bcal}{\mathcal{B}}
\newcommand{\Tcal}{\mathcal{T}}
\newcommand{\Dcal}{\mathcal{D}}

\newcommand{\Scal}{\mathcal{S}}
\newcommand{\sn}{{\mathcal S}^n}
\newcommand{\snp}{{\mathcal S}^n_+}
\newcommand{\snml}{{\mathcal S}^{n-\ell}}
\newcommand{\R}{\mathbb{R}}

\newcommand{\CL}{\textrm{CL}}
\newcommand{\ML}{\textrm{ML}}
\newcommand{\textif}{\textrm{if}}
\newcommand{\st}{\textrm{s.t.}}
\newcommand{\trace}{\textrm{trace}}
\newcommand{\Diag}{\textrm{Diag}}

\begin{document}

\ifJOC
\else 

\title{SOS-SDP: an Exact Solver for Minimum Sum-of-Squares Clustering}
\date{\today}
\author{Veronica Piccialli, Antonio M.~Sudoso, Angelika Wiegele}

\fi

\ifJOC
\else    
\maketitle
\begin{abstract}
The minimum sum-of-squares clustering problem (MSSC) consists of partitioning $n$ observations into $k$ clusters in order to minimize the sum of squared distances from the points to the centroid of their cluster. In this paper, we propose an exact algorithm for the MSSC problem based on the branch-and-bound technique. The lower bound is computed by using a cutting-plane procedure where valid inequalities are iteratively added to the Peng-Wei SDP relaxation. The upper bound is computed with the constrained version of $k$-means where the initial centroids are extracted from the solution of the SDP relaxation. In the branch-and-bound procedure, we incorporate instance-level must-link and cannot-link constraints to express knowledge about which data points should or should not be grouped together. We manage to reduce the size of the problem at each level preserving the structure of the SDP problem itself. The obtained results show that the approach allows to successfully solve for the first time real-world instances up to 4000 data points.
\end{abstract}
\fi 

\section{Introduction}\label{sec:intro}
Clustering is the task of partitioning a set
of objects into homogeneous and/or well-separated groups, called clusters. Cluster analysis is the discipline that studies methods and algorithms for clustering objects according to a suitable similarity measure. It belongs to unsupervised learning since it does not use class labels. Two main clustering approaches exist: hierarchical clustering, which assumes a tree structure in the data and builds nested clusters, and partitional clustering. Partitional clustering generates all the clusters at the same time without assuming a nested structure. Among partitional clustering, the minimum sum-of-squares clustering problem (MSSC) or sum-of-squares (SOS) clustering, is one of the most popular and well studied. MSSC asks to partition $n$ given data points into $k$ clusters so that the sum of the Euclidean distances from each data point to the cluster centroid is minimized. 

The MSSC commonly arises in a wide range of disciplines and applications, as for example image segmentation \citep{dhanachandra2015image, shi2000normalized}, credit risk evaluation \citep{CARUSO2021100850}, biology \citep{jiang2004cluster}, customer segmentation \citep{syakur2018integration}, document clustering \citep{mahdavi2009harmony}, and as a technique for the missing values imputation \citep{zhang2006clustering}. 


The MSSC can be stated as follows for fixed $k$:
\begin{subequations}
\label{eq:MSSC}
\begin{align}
\min~ & \sum_{i=1}^n \sum_{j=1}^k x_{ij}\|p_i - c_j\|^2  \\
\st~ & \sum_{j=1}^k x_{ij} = 1,\quad \forall i \in \{1,\dots,n\} \label{eq:MSSCa}\\
& \sum_{i=1}^n x_{ij} \ge 1, \quad \forall j \in \{1,\dots,k\}\label{eq:MSSCb}\\
& x_{ij} \in \{0,1\}, \quad \forall i \in \{1,\dots,n\}\; \forall j \in \{1,\dots,k\} \\
& c_j \in \mathbb{R}^d, \quad \forall j \in \{1, \dots, k \}.
\end{align}
\end{subequations}
Here, $p_i \in \R^d$, where $d$ is the number of features, $i \in \{1,\dots,n\}$, are the data points, and the centers of the $k$ clusters are at the (unknown) points $c_j$, $j\in \{1,\dots,k\}$. 
For convenience, we sometimes collect all the data points $p_i$ as rows in a matrix $W_p$.
The binary decision variable $x_{ij}$ expresses whether data point $i$ is assigned to cluster $j$ or not.
Constraints~\eqref{eq:MSSCa} make sure that each point is assigned to a cluster, and constraints~\eqref{eq:MSSCb} guarantee that none of the $k$ clusters is empty.

Setting the gradient of the objective function with respect to $c$ to zero yields
\begin{equation*}
\sum_{i=1}^n x_{ij}(c^r_j - p^r_i) = 0,\quad \forall j \in \{1,\dots,k\}\; \forall r \in \{1,\dots,d\}
\end{equation*}
and we obtain the formula for the point in the center of each cluster 
\begin{equation*}
c^r_j = \frac{\sum_{i=1}^n x_{ij} p^r_i}{\sum_{i=1}^n x_{ij}}, \quad \forall j \in \{1,\dots,k\}\; \forall r \in \{1,\dots,d\}.
\end{equation*}
Replacing the formula for $c$ in~\eqref{eq:MSSC}, we get 
\begin{subequations}
\label{eq:MSSC2}
\begin{align}
\min~ & \sum_{i=1}^n \sum_{j=1}^k x_{ij}\Big\|p_i - \frac{\sum_{l=1}^n x_{lj} p_l}{\sum_{l=1}^n x_{lj}}\Big\|^2  \\
\st~ & \sum_{j=1}^k x_{ij} = 1,\quad \forall i \in \{1,\dots,n\}\\
& \sum_{i=1}^n x_{ij} \ge 1, \quad \forall j \in \{1,\dots,k\}\\
& x_{ij} \in \{0,1\}, \quad \forall i \in \{1,\dots,n\}\; \forall j \in \{1,\dots,k\}.
\end{align}
\end{subequations}

\subsection{Literature Review}
The MSSC is known to be NP-hard in $\mathbb{R}^2$ for general values of $k$ \citep{mahajan2012planar}, and in higher dimension even for $k=2$ \citep{aloise2009np}. The one-dimensional case is proven to be solvable in polynomial time. In particular, \cite{wang2011ckmeans} proposed an $O(kn^2)$ time and $O(kn)$ space dynamic programming algorithm for solving this special case. Because of MSSC's computational complexity, heuristic approaches and approximate algorithms are usually preferred over exact methods. 

The most popular heuristic for solving MSSC is $k$-means \citep{macqueen1967some, lloyd1982least},
that alternates the centroid initialization with the assignments of points until centroids do not move anymore. The main disadvantage of $k$-means is that it  produces locally optimal solutions that can be far from the global minimum, and it is extremely sensitive to the initial assignment of centroids. For this reason, a lot of research has been dedicated to finding efficient initialization for $k$-means (see for example \cite{arthur2006k,improvedkmeans2018,franti2019much} and references therein). However, an efficient initialization may not be enough in some instances, so that different strategies have been implemented in order to improve the exploration capability of the algorithm. A variety of heuristics and metaheuristics have been proposed, following the standard metaheuristic framework, e.g., simulated annealing \citep{lee2021simulated}, tabu search \citep{ALSULTAN19951443}, variable neighborhood search \citep{HANSEN2001405,Orlov2018}, iterated local search \citep{likas2003global}, evolutionary algorithms \citep{MAULIK20001455,SARKAR1997975}). In the work of \cite{tao2014new,BAGIROV201612,KARMITSA2017367,KARMITSA2018245},  DC (Difference of Convex functions) programming is used to define efficient heuristic algorithms for clustering large datasets. 
The algorithm $k$-means has also been used as a local search subroutine in different algorithms, as in the population-based metaheuristic developed in \cite{gribel2019hg} and in the differential evolution scheme proposed in \cite{Schoen2021}.

Recently, thanks to the enhancements in computers' computational power and to the progress in mathematical programming, the exact resolution of MSSC has become way more achievable. In this direction, mathematical programming algorithms based on branch-and-bound and column generation have produced guaranteed globally optimal solutions for small and medium scale instances. Due to the NP-hardness of the MSSC, the computational time of globally optimal algorithms quickly increases  with the size of the problem. However, besides the importance of finding optimal solutions for some clustering applications, certified optimal solutions remain extremely valuable as a benchmark tool since they can be used for evaluating, improving, and developing heuristics and approximate methods.
Compared to the huge number of papers proposing heuristics and approximate methods for the MSSC problem, the number of articles proposing exact algorithms is much smaller.  

One of the earliest attempts was the integer programming formulation proposed by \citet{rao1971cluster}, which requires the cluster sizes to be fixed in advance and is limited to small instances. A first branch-and-bound algorithm was proposed by \citet{koontz1975branch} and extended by \citet{diehr1985evaluation}.
The idea is to use partial clustering solutions on a subset $S$ of the main dataset $D$ to determine improved bounds and clusters on the entire sample by a branch-and-bound search. The key observation is that the optimal objective function value of the MSSC on $D$ is greater or equal than the optimal objective function value of the MSSC on $S$ plus the optimal objective function value of the MSSC on $D - S$.
This approach was later improved by \citet{brusco2006repetitive},
who developed a repetitive-branch-and-bound algorithm (RBBA). After a proper reordering of the entities in $D$, RBBA  solves a sequence of subproblems of increasing size with the branch-and-bound technique. While performing a branch-and-bound for a certain subproblem, Brusco's algorithm exploits the optimal solutions found for the previous subproblems which provide tighter bounds compared to the ones used by \cite{koontz1975branch} and \cite{diehr1985evaluation}.  RBBA provided optimal solutions for well separated synthetic datasets with up to 240 objects. Poorly separated problems with no inherent cluster structure were optimally solved for up to 60 objects. 
\citet{sherali2005global} proposed a different branch-and-bound algorithm where tight lower bounds are determined by using the reformulation-linearization-technique (RLT), see \citet{sherali1998reformulation}. The authors claim that this algorithm allows for the exact resolution of problems of size up to 1000 entities, but those results seem to be hard to reproduce. The computing times in an attempted replication by \citet{aloise2011evaluating} were already high for real datasets with about 20 objects.

A column generation algorithm for MSSC was proposed by \citet{du1999interior}. The master problem is solved by an interior point method, whereas the auxiliary problem of finding a column with negative reduced cost is expressed as a hyperbolic program with binary variables. Variable-neighborhood-search heuristics are used to find a good initial solution 
and to accelerate the resolution of the auxiliary problem. This approach has been considered a successful one, since it solved for the first time medium size benchmark instances (i.e., instances with 100--200 entities), including the popular Iris dataset, which encounters 150 entities. However, the bottleneck of the algorithm lies in the resolution of the auxiliary problem, and more precisely, in the unconstrained quadratic 0-1 optimization problem. Later this algorithm was further improved by  \citet{aloise2012improved} who define a different geometric-based approach for solving the auxiliary problem. In particular, the solution of the auxiliary problem is achieved by solving a certain number of convex quadratic problems. If the points to be clustered are in the plane, the maximum number of convex problems to solve is polynomially bounded. When the points are not in the plane, in order to solve the auxiliary problems the cliques in a certain graph (induced by the current solution of the master problem) have to be found. The algorithm is more efficient when the graph is sparse, and the graph becomes sparser when the number of clusters $k$ increases. Therefore, the algorithm proposed in \citet{aloise2012improved} is particularly efficient in the plane and when $k$ is large. Their method was able to provide exact solutions for large scale problems, including one instance of 2300 entities when the ratio between $n$ and $k$ is small. 

Recently, \citet{peng2007approximating} by using matrix arguments proved the equivalence between the MSSC formulation and a model called 0-1 semidefinite programming (SDP), in which the eigenvalues of the matrix variable are binary. Using this result, \citet{aloise2009branch} proposed a branch-and-cut algorithm for MSSC where lower bounds are obtained from the linear programming relaxation of the 0-1 SDP model. This algorithm manages to obtain exact solutions for datasets up to 200 entities with computing times comparable with those obtained by the column generation method proposed by \citet{du1999interior}.

Constant-factor approximation algorithms have also been developed in the literature, both for fixed number of clusters $k$ and for fixed dimension $d$ \citep{kanungo2004local}. Among these methods, \citet{peng2007approximating} proposed a rounding procedure to extract a feasible solution of the original MSSC 
from the approximate solution of the relaxed SDP problem. More in detail, they use
the Principal Component Analysis (PCA) to reduce the dimension of the dataset and then perform clustering on the projected PCA space. They showed that this algorithm can provide a 2-approximate solution to the MSSC. 
More recently, \citet{prasad2018improved} proposed a new approximation algorithm that utilizes an improved copositive conic reformulation of the MSSC. Starting from this reformulation, the authors derived a hierarchy of accurate SDP relaxations obtained by replacing the completely positive cone with progressively tighter semidefinite outer approximations. Their SDP relaxations provide better lower bounds than the Peng-Wei one but do not scale well when the size of the problem increases. 

\subsection*{Main results and outline}
The main contributions of this paper are the following:
\begin{description}
    \item[(i)] we define the first SDP based branch-and-bound algorithm for MSSC, and we use a cutting-plane procedure for strengthening the bound, following a recent strand of research \citep{demeijer2021sdpbased};
    \item[(ii)] we define a shrinking procedure that allows reducing the size of the problem when introducing must link constraints;
    \item[(iii)] we exploit the SDP solution for a smart initialization of the constrained version of $k$-means that yields high quality upper bounds;
    \item[(iv)] for the first time, we manage to find the exact solution for instances of size up to $n=4000$.
\end{description}

This paper is structured as follows. In Section~\ref{sec:bound} we introduce equivalent formulations for the MSSC and derive relaxations based on semidefinite programming (SDP). 
In Section~\ref{sec:branching} we analyze the SDP problems that arise at each node within the branch-and-bound tree and discuss the selection of the branching variable.
In Section~\ref{sec:bab} the details about the bound computation are discussed, including a post-processing procedure that produces a ``safe'' bound from an SDP that is solved to medium precision only.
Section~\ref{sec:heuristic} gives all the details on the heuristic used to generate feasible clusterings. 
The details of our implementation and exhaustive numerical results are presented in  Section~\ref{sec:numericalresults}.
Finally, Section~\ref{sec:conclusion} concludes the paper.

\subsection*{Notation}
Let $\sn$ denote the set of all $n\times n$ real symmetric matrices. We denote by $M\succeq 0$ that matrix $M$ is positive semidefinite and let ${\mathcal S}_+^n$ be the set of all positive semidefinite matrices of order
$n\times n$. We denote by $\inprod{\cdot}{\cdot}$ the 
trace inner product. That is, for any
$M, N \in \mathbf{R}^{n\times n}$, we define $\inprod{M}{N}:= \trace (M^\top N )$. Its associated norm is the Frobenius norm, denoted by $\| M\|_F := \sqrt{\trace (M^\top M )}$. We define the linear map $\Acal: \sn \rightarrow \mathbb{R}^{m_1}$ as $(\Acal(X))_i = \inprod{A_i}{X}$, where $A_i \in \sn$, $i=1,\dots,m_1$, and the linear map $\Bcal: \sn \rightarrow \mathbb{R}^{m_2}$ as $(\Bcal(X))_i = \inprod{B_i}{X}$, where $B_i \in \sn$, $i=1,\dots,m_2$.

We define by $e_n$ the vector of all ones of length $n$. We omit the subscript in case the dimension is clear from the context.
We denote by $E_{i}$ the symmetric matrix such that $\inprod{E_{i}}{Z}$ is the sum of row~$i$ of $Z$.

\section{A Lower Bound based on Semidefinite Programming}\label{sec:bound}
We briefly remind the Peng-Wei SDP relaxation to Problem~\eqref{eq:MSSC2} that will be the basis of the bounding procedures within our exact algorithm.

Consider matrix $W$ where the entries are the inner products of the data points, i.e., $W_{ij} = p_i^\top p_j$ for $i,j \in \{1,\dots,n\}$. Furthermore, collect the binary decision variables $x_{ij}$ from~\eqref{eq:MSSC2} in the $n\times k$ matrix $X$ and define matrix $Z$ as
\begin{equation*}
Z = X(X^\top X)^{-1}X^\top.
\end{equation*}
\citet{peng2007approximating} introduced a different but equivalent formulation for the MSSC, yielding the following
optimization problem:
\begin{subequations}
\label{eq:PengSDP}
\begin{align}
\min~ & \inprod{-W}{Z}  \\
\st~ & Ze = e\\
& \textrm{tr}(Z) = k\\
& Z \ge 0, \ Z^2 = Z, \ Z = Z^\top.
\end{align}
\end{subequations}
We can convert Problem~\eqref{eq:PengSDP} into a rank constrained optimization problem. In fact we can replace the constraints $Z^2 = Z$ and $Z = Z^\top$ with a rank constraint and a positive semidefiniteness constraint on $Z$, yielding the following problem:
\begin{subequations}
\label{eq:RankSDP}
\begin{align}
\min~ & \inprod{-W}{Z}  \\
\st~ & Ze = e\\
& \textrm{tr}(Z) = k\\
& Z \ge 0, \ Z \in \snp\\
& \textrm{rank}(Z) = k.
\end{align}
\end{subequations}
In order to prove the equivalence of Problems~\eqref{eq:PengSDP} and~\eqref{eq:RankSDP}, we need the definition of an idempotent matrix and its characterization in terms of eigenvalues given by Lemma~\ref{lemma:ideig}.
\begin{definition}
A symmetric matrix $Z$ is idempotent if $Z^2 = ZZ = Z$.
\end{definition}

\begin{lemma}\label{lemma:ideig}
A symmetric matrix $Z$ is idempotent if and only if all its eigenvalues are either 0 or 1.
\end{lemma}

\ifJOC
\proof{Proof.}
\else
\begin{proof}
\fi
Let $Z$ be idempotent, $\lambda$ be an eigenvalue and $v$ a
corresponding eigenvector then $\lambda v = Zv = ZZv = \lambda Zv = \lambda^2 v$.
Since $v \neq 0$ we find $\lambda - \lambda^2 = \lambda (1 - \lambda) = 0$ so either $\lambda = 0$ or $\lambda = 1$.

To prove the other direction, consider the eigenvalue decomposition of $Z$,
$Z = P \Lambda P^\top$, 
where $\Lambda$ is a diagonal matrix having the eigenvalues $0$ and $1$ on the diagonal, and $P$ is orthogonal. Then, since $\Lambda^2 = \Lambda$, we get
\begin{equation*}
    Z^2 = P \Lambda P^\top P \Lambda P^\top = P \Lambda^2 P^\top = P \Lambda P^\top = Z.
\end{equation*}
\ifJOC
\hfill\Halmos
\endproof
\else
\end{proof}
\fi


\begin{theorem}
Problems~\eqref{eq:PengSDP} and~\eqref{eq:RankSDP} are equivalent.
\end{theorem}

\ifJOC
\proof{Proof.}
\else
\begin{proof}
\fi
Let $Z$ be a feasible solution of Problem \eqref{eq:PengSDP}. We first show that $Z^2 = Z$ and $Z = Z^T$ imply $Z \in \snp$. In fact, for all $v$ we have:
\begin{equation*}
    v^\top Z v = v^\top Z^2 v = v^\top Z Z v = v^\top Z (v^\top Z^\top)^\top = (v^\top Z) (v^\top Z)^\top = \| v^\top Z \|_2^2 \geq 0.
\end{equation*}
Since $Z$ is symmetric idempotent, the number of eigenvalues equal to 1 is $\textrm{tr}(Z) = \textrm{rank}(Z) = k$.

To prove the other direction, let $Z$ be a feasible solution of Problem \eqref{eq:RankSDP}. If $\textrm{rank}(Z) = k$, then $Z$ has $n-k$ eigenvalues equal to 0. Furthermore, let $\lambda_1 \geq \lambda_2 \geq \ldots > \lambda_n\ge 0$ be the eigenvalues of Z, then 

\begin{equation*}
    \textrm{tr}(Z) = \sum_{i=1}^{n} \lambda_i = \sum_{i=1}^{k} \lambda_i + \sum_{i=k+1}^{n} \lambda_i = \sum_{i=1}^{k} \lambda_i = k.
\end{equation*}
Constraints $Z\succeq 0$, $Z\ge 0$ and $Ze=e$ imply that the eigenvalues of $Z$ are bounded by one (see, e.g., Lemma~\ref{lem:eigboundZ}). Hence, the trace constraint is satisfied if and only if the positive eigenvalues are all equal to~1. This shows that $\lambda(Z) \in \{0, 1\}$ and therefore $Z$ is symmetric idempotent.
\ifJOC
\par \hfill\Halmos
\endproof
\else
\end{proof}
\fi

By dropping the non-convex rank constraint from Problem~\eqref{eq:RankSDP}, we obtain the SDP relaxation which is the convex optimization problem
\begin{subequations}
\label{eq:SDP}
\begin{align}
\min~ & \inprod{-W}{Z}  \\
\st~ & Ze = e\\
& \textrm{tr}(Z) = k\\
& Z \ge 0, \ Z \in \snp
\end{align}
\end{subequations}

\subsection{Strengthening the Bound through Inequalities}

The SDP relaxation~\eqref{eq:SDP} can be tightened by adding valid inequalities and solving the resulting SDP in a cutting-plane fashion. In this section, we present the class of inequalities we use for strengthening the bound. For each class, we describe the separation routine used.

We consider three different sets of inequalities:

\begin{description}
    \item[Pair inequalities.] In any feasible solution of~\eqref{eq:RankSDP}, it holds that \begin{equation}\label{eq:pairs}
        Z_{ij}\le Z_{ii},\quad Z_{ij}\le Z_{jj}\quad \forall i,j \in \{1,\dots,n\}, i\not=j.
    \end{equation} 
    This set of $n(n-1)$ inequalities were used by \citet{peng2005new} and in the branch-and-cut proposed by \citet{aloise2009branch}.
    \item[Triangle Inequalities.] 
    The triangle inequalities are based on the observation that if points $i$ and $j$ are in the same cluster and points $j$ and $h$ are in the same cluster, then points $i$ and $h$ necessarily must be in the same cluster. The resulting $3\binom{n}{3}$ inequalities are:
    \begin{equation}\label{eq:triangle}
            Z_{ij}+Z_{ih}\le Z_{ii}+Z_{jh}\quad \forall i,j,h \in \{1,\dots,n\}, i,j,h ~\mathrm{distinct}.
    \end{equation}
    
    These inequalities were already introduced by \citet{peng2005new}, and used also by \citet{aloise2009branch}. 
    \item[Clique Inequalities.] If the number of clusters is $k$, for any subset $Q$ of $k+1$ points at least two points have to be in the same cluster (meaning that at least one $Z_{ij}$ needs to be positive and equal to $Z_{ii}$ for all $(i,j)\in Q$). This can be enforced by the following inequalities:
       \begin{equation}\label{eq:clique}
        \sum_{(i,j)\in Q,i<j}Z_{ij}\ge \frac{1}{n-k+1} \quad\forall Q\subset\{1,\ldots,n\},\,|Q|=k+1.
    \end{equation}
    These $\binom{n}{k+1}$ inequalities are similar to the clique inequalities for the $k$-partitioning problem \citep{chopra1993partition}, the difference lies in the right hand side, that in that case is equal to~1, whereas here we use the smallest possible value that an element on the diagonal of $Z$ can hold. 
    \end{description}
    

Pair and triangle inequalities are known to be valid for Problem~\eqref{eq:PengSDP}, see \cite{peng2005new} and \cite{de2020ratio}. 
It remains to show that also the clique inequalities are valid.

\begin{lemma}
The clique inequalities~\eqref{eq:clique} are valid for Problem~\eqref{eq:PengSDP}.
\end{lemma}

\ifJOC
\proof{Proof.}
\else
\begin{proof}
\fi
  The left hand side of \eqref{eq:clique} has $\binom{k+1}{2}$ terms, and we know that $Z_{ii}\ge\frac{1}{n-k+1}$, since the cardinality of a cluster can be at most $n-k+1$.
  Given that the number of clusters is $k$, for any set of $k+1$ points at least two points have to be in the same cluster, say points $i$ and $j$.  Then, for any feasible clustering $Z$, at least the element $Z_{ij}$ in the left hand side of \eqref{eq:clique} needs to be different from zero, therefore equal to $Z_{ii}$, and hence \eqref{eq:clique} must hold. 
\ifJOC
\hfill\Halmos
\endproof
\else
\end{proof}
\fi

\section{Branching: Subproblems within a Branch-and-Bound Algorithm and Variable Selection}\label{sec:branching}
Our final goal is to develop a branch-and-bound scheme to solve the MSSC to optimality using relaxation~\eqref{eq:SDP} strengthened by some of the inequalities~\eqref{eq:pairs}--\eqref{eq:clique}. In this section we examine the problems that arise after branching. To keep the presentation simple and since everything carries over in a straightforward way, we omit in this section the inclusion of inequalities~\eqref{eq:pairs}--\eqref{eq:clique}.

The branching decisions are as follows. Given a pair $(i,j)$, 
\begin{itemize}
    \item points $p_i$ and $p_j$ should be in different clusters, i.e., they  \textit{cannot link} or
    \item points $p_i$ and $p_j$ should be in the same cluster, i.e., they \textit{must link}.
\end{itemize}
By adding constraints due to the branching decisions, the problem changes. However, the structure of the SDP remains similar. 
In this section we describe the subproblems to be solved at each node in the branch-and-bound tree. Each such SDP is of the form
\begin{subequations}
\label{eq:SDPbab}
\begin{align} 
\min~ & \inprod{-\Tcal^{\ell} W (\Tcal^{\ell})^\top}{Z^{\ell}}  \\
\st~ & Z^{\ell} e^{\ell} = e \\
& \inprod{\Diag(e^{\ell})}{Z^{\ell}} = k\\
& Z^\ell_{ij} = 0 \quad (i,j) \in \CL\\
&Z^{\ell} \ge 0, \ Z^{\ell} \in \Scal^+_{n-\ell}
\end{align}
\end{subequations}
where $\CL$ (cannot link) is the set of pairs that must be in different clusters and matrix $\Tcal^{\ell}$ and vector $e^\ell$ encode the branching decisions that ask data points to be in the same cluster (i.e., they must link). We describe this in detail in the subsequent sections.

\subsection{Branching Decisons}

In case we want to have $i$ and $j$ in different clusters, we add the constraint $Z_{ij} = 0$ to the SDP, i.e., we add the pair $(i,j)$ to the set $\CL$.

In the other case, i.e., when the decison is to have $i$ and $j$ in the same cluster, we proceed as follows.
Assume at the current node we have $n$ points and we decide that on this branch the two points $p_i$ and $p_j$ have to be in the same cluster. We can reduce the size of  $W_p$ (the matrix having data points $p_i$ as rows) by substituting row $i$ by $p_{i}+p_j$ and omitting row $j$. To formalize this procedure, we introduce the following notation.

Let $b(r) = (i,j)$, $i<j$, be the branching pair in branching decision at level $r$ and $b(1),\dots,b(\ell)$ a sequence of consecutive branching decisions. Furthermore, let $g(r)=(\underline{i}, \underline{j})$ be the corresponding global indices.

Define $\Tcal^{\ell} \in \{0,1\}^{(n-\ell) \times n}$ as 
\[
\Tcal^{\ell} = T^{b(\ell)} T^{b(\ell-1)} \dots T^{b(1)}
\]
where the $(n-r)\times (n-r+1)$ matrix $T^{b(r)}$ for branching decision $b(r)=(i,j)$ is defined by
\[
T^{b(r)}_{s,\cdot} = \left\{ \begin{array}{ll}
u_s & \textif~ 1 \le s < i ~\textrm{and}~ i+1\le s\le j\\
u_i + u_j & \textif~ s=i\\
u_{s+1} & \textif~ j<s\le n-r
\end{array}\right.
\]
with $u_s$ being the unit vector of size $(n-r+1)$.
Furthermore, we define $T^{b(0)} = I_n$. 

Note that $T^{b(r)}\cdot M$ builds a matrix of size $(n-r)\times (n-r+1)$ by adding rows $i$ and $j$ of $M$ and putting the result into row $i$ while row $j$ is removed and all other rows remain the same.

We also define the vector $e^\ell \in \R^{n-\ell}$ as
\[
e^\ell = \Tcal^\ell e
\]
where $e$ is the vector of all ones of length $n$. 

\begin{remark}\label{rem:ell}
Note that in $(e^\ell)$ the number of points that have been fixed to belong to the same cluster along the branching decisions $b(1),\dots,b(\ell)$ are given. 
Furthermore,  $\Tcal^\ell(\Tcal^{\ell})^\top = \Diag(e^\ell)$.
\end{remark}

We now show that this shrinking operation corresponds to the must-link branching decisions. Consider the following two semidefinite programs.

\begin{subequations}
\label{eq:SDPell}
\begin{align} 
\min~ & -\inprod{\Tcal^{\ell} W (\Tcal^{\ell})^\top}{Z^{\ell}}  \\
\st~ & Z^{\ell} e^{\ell} = e_{n-\ell} \label{eq:SDPellb}\\
& \inprod{\Tcal^{\ell}(\Tcal^{\ell})^\top}{Z^{\ell}} = k \label{eq:SDPellk}\\
&Z^{\ell} \ge 0, Z^{\ell} \in \Scal^+_{n-\ell}
\end{align}
\end{subequations}

and 

\begin{subequations}
\label{eq:SDPbranch}
\begin{align}
\min~ & -\inprod{W}{Z}  \\
\st~ & Ze = e\\
& \inprod{I}{Z} = k\\
& Z_{i\cdot} = Z_{j\cdot} \quad \forall \{i,j\} \in g(l), ~ l \in \{1,\dots, \ell\} \label{eq:SDbranch-rows}\\
&Z \ge 0, Z \in \snp
\end{align}
\end{subequations}

\begin{theorem}
Problems~\eqref{eq:SDPell} and~\eqref{eq:SDPbranch} are equivalent.
\end{theorem}
\ifJOC
\proof{Proof.}
\else
\begin{proof}
\fi
Let $Z^{\ell}$ be a feasible solution of Problem~\eqref{eq:SDPell}. Define $Z = (\Tcal^{\ell})^\top Z^{\ell} \Tcal^{\ell}$. This is equivalent to expanding the matrix by replicating the rows according to branching decisions. Therefore, \eqref{eq:SDbranch-rows} holds by construction.
Clearly, $Z \ge 0$ and $Z\in \snp$ hold as well. Moreover, we have that 

\[\inprod{I}{Z} =  \inprod{I}{(\Tcal^{\ell})^\top Z^{\ell} \Tcal^{\ell}} = \inprod{\Tcal^{\ell}(\Tcal^{\ell})^\top}{Z^{\ell}} = k
\]
and 
\[
Ze = (\Tcal^{\ell})^\top Z^{\ell} \Tcal^{\ell} e = (\Tcal^{\ell})^\top Z^{\ell} e^{\ell} = 
(\Tcal^{\ell})^\top e_{n-\ell} = e_n.
\]
Furthermore, 
\[ \inprod{W}{Z} = \inprod{W}{(\Tcal^{\ell})^\top Z^{\ell} \Tcal^{\ell}}
= \inprod{(\Tcal^{\ell})W(\Tcal^{\ell})^\top}{Z^{\ell}}
\]
and thus $Z$ is a feasible solution of Problem~\eqref{eq:SDPbranch} and the values of the objective functions coincide.

We next prove that any feasible solution of Problem~\eqref{eq:SDPbranch} can be transformed into a feasible solution of Problem~\eqref{eq:SDPell} with the same objective function value. In order to do so, we define the matrix 
\[ \Dcal^\ell = \Diag(1/e^\ell)\]
where $1/e^\ell$ denotes the vector that takes the inverse elementwise. It is straightforward to check that
\[
\Dcal^\ell \Tcal^\ell (\Tcal^\ell)^\top \Dcal^\ell = \Dcal^\ell.
\]
Assume that $Z$ is a feasible solution of Problem~\eqref{eq:SDPbranch}
and set $Z^{\ell} = \Dcal^\ell \Tcal^{\ell} Z (\Tcal^{\ell})^\top \Dcal^\ell$. If $Z$ is nonnegative and positive semidefinite, then so is $Z^\ell$.
Furthermore, we can derive 
\begin{align*} 
\inprod{\Tcal^{\ell}(\Tcal^{\ell})^\top}{Z^{\ell}} & =\inprod{\Tcal^{\ell}(\Tcal^{\ell})^\top}{\Dcal^\ell \Tcal^{\ell} Z (\Tcal^\ell)^\top \Dcal^\ell}\\
&= \inprod{\Dcal^\ell\Tcal^{\ell}(\Tcal^{\ell})^\top \Dcal^\ell}{\Tcal^\ell Z (\Tcal^\ell)^\top }\\
&= \inprod{\Dcal^\ell}{ \Tcal^\ell Z (\Tcal^\ell)^\top }
= \sum_{l=1}^{n-\ell} \frac{1}{e^\ell_l} \sum_{j \in g(l)} \sum_{i\in g(l)}  Z_{ij}\\
&\stackrel{(*)}{=} \sum_{l=1}^{n-\ell} \frac{1}{e^\ell_l} \sum_{j \in g(l)} \sum_{i\in g(l)}  Z_{ii}
= \sum_{l=1}^{n-\ell} \frac{1}{e^\ell_l} e^\ell_l \sum_{i\in g(l)}  Z_{ii}\\
&= \sum_{l=1}^{n-\ell} \sum_{i\in g(l)}  Z_{ii}= \sum_{i=1}^n  Z_{ii}= k.
\end{align*}
Note that the equality~$(*)$ holds since $Z_{i,j} = Z_{r,s}$ for any $i,j,r,s \in g(l)$. This ensures that constraint~\eqref{eq:SDPellk} holds for $Z^\ell$.

To prove~\eqref{eq:SDPellb} consider the equations 
\begin{align*}
Z^\ell e^\ell & = \Dcal^\ell \Tcal^{\ell} Z (\Tcal^{\ell})^\top \Dcal^\ell e^\ell
= \Dcal^\ell \Tcal^{\ell} Z (\Tcal^{\ell})^\top e_{n-\ell} \\
&= \Dcal^\ell \Tcal^{\ell} Z e
= \Dcal^\ell \Tcal^{\ell} e = \Dcal^\ell e^\ell = e_{n-\ell}.
\end{align*}

It remains to show that the objective function values coincide.
\begin{align*}
\inprod{\Tcal^{\ell} W (\Tcal^{\ell})^\top}{Z^{\ell}} 
&= 
\inprod{\Tcal^{\ell} W (\Tcal^{\ell})^\top}{\Dcal^\ell \Tcal^{\ell} Z (\Tcal^{\ell})^\top \Dcal^\ell}\\
&= 
\inprod{ W }{(\Tcal^{\ell})^\top\Dcal^\ell \Tcal^{\ell} Z (\Tcal^{\ell})^\top \Dcal^\ell\Tcal^{\ell}}\\
&= 
\inprod{ W }{Z}.\\
\end{align*}
As for the last equation, note that pre- and postmultiplying $Z$ by $(\Tcal^{\ell})^\top\Dcal^\ell \Tcal^{\ell}$ ``averages'' over the respective rows of matrix $Z$. Since these respective rows are identical due to~\eqref{eq:SDbranch-rows}, the last equation holds.
\ifJOC
\hfill\Halmos
\endproof
\else
\end{proof}
\fi

\begin{remark}
  The addition of constraints $Z_{ij}=0$ for datapoints $i,j$ that should not belong to the same cluster also goes through in the above equivalence. However, to keep the presentation simple we did not include it in the statement of the theorem above.
\end{remark}

\begin{remark}
It is straightforward to include the additional constraints~\eqref{eq:pairs}, \eqref{eq:triangle}, and~\eqref{eq:clique} in the subproblems, i.e., in case of shrinking the problem, the constraints are still valid. Again, to keep notation simple, we omitted these constraints in the presentation above.
Further discussions on including these inequalities are in Section~\ref{sec:boundcomp}.
\end{remark}

\subsection{Variable Selection for Branching}\label{sec:variableselection}
In a matrix $Z$ corresponding to a clustering, for each pair $(i,j)$ either $Z_{ij}=0$ or $Z_{ii} = Z_{ij}$. \citet{peng2005new} propose a simple branching scheme. Suppose that for the optimal solution of the SDP relaxation there are indices $i$ and $j$ such that $Z_{ij}(Z_{ii}-Z_{ij}) \neq 0$ then one can produce a cannot-link branch with $Z_{ij} = 0$ and a must-link branch with $Z_{ii} = Z_{ij}$. Regarding the variable selection the idea is to choose indices $i$ and $j$ such that in both branches we expect a significant improvement  of the lower bound. In~\cite{peng2005new} the branching pair is chosen as the  \[\argmax_{i,j} \{ \min \{Z_{ij},Z_{ii}-Z_{ij}\} \}.\]   
Here we propose a variable selection strategy that is coherent with the way we generate the cannot-link and the must-link subproblems. In fact, we observe that in a matrix $Z$ corresponding to a clustering, for each pair $(i, j)$ either $Z_{ij} = 0$ or $Z_{i\cdot} = Z_{j\cdot}$. This motivates the following strategy to select a pair of data points to branch on
\[\argmax_{i,j} \{ \min \{Z_{ij}, \|Z_{i\cdot} - Z_{j\cdot}\|_2^2 \} \}.\]   

In case this maximizer gives a value close to zero, say $10^{-5}$, the SDP solution corresponds to a feasible clustering. 

\subsubsection*{Variable selection on the shrunk problem}

The strategy for the variable selection still carries over on the shrunk problems.
Since $Z$ is obtained from $Z^\ell$ only by repeating rows and columns, every pair $(Z^\ell_{ij}, Z^\ell_{ii})$ appears also in $Z$ and vice versa.
Moreover, within the already merged points, by construction $Z^\ell_{ii}=Z^\ell_{ij}$ and hence this can never be a branching candidate again.


\section{Branch-and-Bound Algorithm}\label{sec:bab}

We now put the bound computation (see Sections~\ref{sec:bound} and~\ref{sec:branching}) together with our way of branching (see Section~\ref{sec:variableselection}) to form our algorithm \texttt{SOS-SDP}. The final ingredient, a heuristic for providing upper bounds, is described in Section~\ref{sec:heuristic}.

\subsection{The Bound Computation}\label{sec:boundcomp}
In order to obtain a strong lower bound, we solve the SDP relaxation~\eqref{eq:SDP} strengthened by the inequalities given in Section~\ref{sec:bound}. 

The enumeration of all pair and triangle inequalities is computationally intractable even for medium size instances. Therefore we use a similar separation routine for both types of inequalities:
\begin{enumerate}
    \item Generate randomly up to $t$ inequalities violated by at least $\varepsilon_{\mathrm{viol}}$
    \item Sort the  $t$ inequalities by decreasing violation
    \item Add to the current bounding problem the $p\ll t$ most violated ones.
\end{enumerate}

As for the clique inequalities, we use the heuristic separation routine described in \cite{ghaddar2011branch} for the minimum $k$-partition problem, that returns at most $n$ valid clique inequalities. More in detail, at each cutting-plane iteration, these cuts are determined by finding $n$ subsets $Q$ with a greedy principle. For each point $i \in S = \{1, \dots, n\}$, $Q$ is initialized as $Q = \{i\}$. Then, until the cardinality of $Q$ does not reach the size $k+1$, $Q$ is updated as $Q = Q \cup \{\argmin_{j \in S \setminus Q} \sum_{q \in Q} Z_{qj} \}$.

We denote by $\Acal(Z^\ell) = b$ the equations from the must-link and cannot-link constraints and by $l \leq \Bcal(Z^\ell) \leq u$ the inequalities representing the cutting planes. 
The cutting-plane procedure performed at each node is outlined in  Algorithm~\ref{alg:cpproc}. 

We stop the procedure when we reach the maximum number of iterations $cp_\textrm{max}$. Another stopping criterion is based on the relative variation of the bound between two consecutive iterations. If the variation is lower than a tolerance $\varepsilon_{\textrm{cp}}$, the cutting-pane method terminates, and we branch.

At each node, we use a cuts inheritance procedure to quickly retrieve several effective inequalities from the parent node and save a significant number of cutting-plane iterations during the bound computation of the children. More in detail, the inequalities that were included in the parent node during the last cutting-plane iteration are passed to its children and included in their problem from the beginning. 
While inheriting inequalities in the $(i,j)$ must link child, the shrinking procedure must be taken into account, updating the indices in the inequalities involved and deleting inequalities involving both points $i$ and $j$. 

In addition to the cuts inheritance, we use a cuts management procedure. A standard cutting-plane algorithm expects the valid inequalities not to be touched after having been included. The efficiency of state-of-the-art SDP solvers considerably deteriorates as we add these cuts, especially when solving large scale instances in terms of $n$. For this reason, after solving the current SDP, we remove the constraints that are not active at the optimum. Of course, inactive constraints may become active again in the subsequent cutting-plane iteration, and this operation could prevent the lower bounds from increasing monotonically; however, empirical results show that this situation happens rarely, and in this case, we decide to stop the cutting-plane procedure and we branch. From the practical standpoint, we notice that removing inactive constraints makes a huge difference since it keeps the SDP problem to a computationally tractable size. The result is that each cutting-plane iteration is more lightweight in comparison to the standard version, and this significantly impacts the overall efficiency of our branch-and-bound algorithm. Our strategy turns out to be more efficient than adding cuts only at the root node and inheriting them in the children. Indeed, if we add cuts only at the root node, the number of nodes in the tree increases since the bound does not improve as much as by repeating the separation routine in each node. Even though the single node is faster since only one SDP is solved, the overall computational time increases.

\begin{algorithm}
\SetKw{Init}{Initialization:}
\SetKw{Or}{or}
\SetKw{Stop}{stop;}
\KwData{A subproblem defined through the current set of equalities $\Acal(Z^\ell) = b$, and inequalities $l \leq \Bcal(Z^\ell) \leq u$, the current global upper bound $\varphi$, the maximum number of cutting-plane iterations $cp_{\max}$, the cutting-plane tolerance $\varepsilon_{\mathrm{cp}}$, the cuts violation tolerance $\varepsilon_{\mathrm{viol}}$, and the cuts removal tolerance $\varepsilon_{\mathrm{act}}$.}
\KwResult{A lower bound $\hat{\delta}^\ell$ on the optimal value of the subproblem}
\Init $i \leftarrow 1$, $\delta_0 \leftarrow -\infty$\

\Repeat{no violated inequalities found}{
solve the current SDP relaxation:


\begin{equation*} 
\hat{\delta}_i^\ell = \min \big\{ \inprod{-\Tcal^{\ell} W (\Tcal^{\ell})^\top}{Z^{\ell}} \colon \Acal(Z^\ell) = b, \ l \leq \Bcal(Z^\ell) \leq u, \ Z^\ell \ge 0, \ Z^\ell \in \snml_+ \big\}
\end{equation*}

and let $\hat{Z}_i^\ell$ be the optimizer\;

\If{$\hat{\delta}_i^\ell \geq \varphi$}{
\Stop the node can be pruned\;
}


\If{$i \geq cp_{\max}$ \Or $\frac{| \hat{\delta}_i^\ell - \hat{\delta}_{i-1}^\ell |}{\hat{\delta}_{i-1}} \leq \varepsilon_{\mathrm{cp}}$}{\Stop return the lower bound $\hat{\delta}_i^\ell$ and branch\;}

remove inactive inequalities with tolerance $\varepsilon_{\mathrm{act}}$ by updating $(\Bcal(\cdot), l, u)$\;

apply the separation routines for pair, triangle and clique inequalities with tolerance $\varepsilon_{\mathrm{viol}}$ and add them to $(\Bcal(\cdot), l, u)$\;

\eIf{no violated inequalities found}{\Stop return the lower bound $\hat{\delta}_i^\ell$ and branch\;}
{add the inequalities by updating $(\Bcal(\cdot), l, u)$\; 
set $i \leftarrow i + 1$\;}
}
\caption{The node processing loop in the branch-and-cut algorithm}
\label{alg:cpproc}
\end{algorithm}

\subsection{Post-processing Using Error Bounds}

Using the optimal solution of the SDP relaxation whithin a branch-and-bound framework requires the computation of ``safe'' bounds. Such safe bounds are obtained by solving the SDP to high precision, which, however, is out of reach when using first-order methods.
In order to obtain a safe bound, we run a post-processing procedure where 
we use a method to obtain rigorous lower bounds on the 
optimal value of our SDP relaxation introduced by \citet{JaChayKeil2007}.
Before describing our post-processing, we state a result bounding the eigenvalues of any feasible solution of~\eqref{eq:SDP}.
\begin{lemma}\label{lem:eigboundZ}
Let $Z\succeq 0$ and $Z\ge 0$. Furthermore, let $Ze=e$. Then the eigenvalues of $Z$ are bounded by one.
\end{lemma}
\ifJOC
\proof{Proof.}
\else
\begin{proof}
\fi
Let $\lambda$ be an eigenvalue of $Z$ with eigenvector $v$, i.e., $Zv = \lambda v$. This implies 
\begin{equation*}
     \lambda |v_i| = | \sum_{j=1}^n z_{ij}v_j|  \le \max_{1\le j\le n} |v_j| \sum_{i=1}^n z_{ij} = \max_{1\le j\le n} |v_j| \mbox{ for all } i\in \{1,\dots,n\} 
\end{equation*}
by nonnegativity of $Z$ and since the row sums of $Z$ are one. 
Therefore, the inequality 
\begin{equation*}
    \lambda \le \frac{\max_{1\le j\le n} |v_j|}{|v_i|}
\end{equation*}
holds for all $i\in \{1,\dots,n\}$, and in particular for $i \in \argmax_{1\le j\le n} |v_j|$ which proves $\lambda \le 1$.
\ifJOC
\hfill\Halmos
\endproof
\else
\end{proof}
\fi

We now restate Lemma~3.1 from~\cite{JaChayKeil2007} in our context.

\begin{lemma}\label{lem:jansson}
Let $S$, $Z$ be symmetric matrices that satisfy
$0 \leq \lambda_{\min}(Z)$ and $\lambda_{\max}(Z) \leq \bar{z}$
for some $\bar{z} \in \R$.
Then the inequality
\begin{equation*}
    \left\langle S,Z\right\rangle \geq \bar{z}\sum_{i \colon  \lambda_i(S) <0}\lambda_i(S)
\end{equation*}
holds.
\end{lemma}
\ifJOC
\proof{Proof.}
\else
\begin{proof}
\fi
 Let $S$ have the eigenvalue decomposition 
 $S=Q\Lambda Q^\top$ where $QQ^\top=I$ and $\Lambda=\Diag(\lambda(S))$.
 Then
 \begin{equation*}
      \left\langle S,Z\right\rangle = \inprod{Q\Lambda Q^\top}{Z}  = \inprod{\Lambda}{Q^\top Z Q}
      = \sum_{i=1}^n \lambda_i(S)Q_{\cdot,i}^\top Z Q_{\cdot,i}
 \end{equation*}
 where $Q_{\cdot,i}$ is column $i$ of matrix $Q$.
 Because of the bounds on the eigenvalues of $Z$, we have $0 \leq Q_{\cdot,i}^\top Z Q_{\cdot,i} \leq \bar{z}.$ Therefore
 $\left\langle S,Z\right\rangle \geq \bar{z}\sum_{i \colon \lambda_i <0}\lambda_i(S)$.
\ifJOC
\hfill\Halmos
\endproof
\else
\end{proof}
\fi

\begin{theorem}\label{thm:errorbound}
  Consider the SDP~\eqref{eq:SDP} together with equations $\Acal(Z)=b$ (e.g., from cannot-link constraints) and inequalities $l \le \Bcal(Z) \le u$ (representing cutting planes) 
with optimal objective function value $p^*$. 
Denote the dual variables by $(\tilde{y},\tilde{u},\tilde{v},\tilde{w},\tilde{P})$, with $\tilde{y}\in \R^{n+1}$, $\tilde{u}$, $\tilde{v}$, $\tilde{w}$ being vectors of appropriate size, $\tilde{P}\in \sn$, $\tilde{P} \ge 0$ and set $\tilde{S} = -W - \sum_{i=1}^n \tilde{y}_iE_i - \tilde{y}_{n+1}I - \Acal^\top(\tilde{u}) + \Bcal^\top(\tilde{v}) - \Bcal^\top(\tilde{w}) - \tilde{P}$. Then
\begin{equation*}
    p^* \ge \sum_{i=1}^n\tilde{y}_i + k\tilde{y}_{n+1} + b^\top \tilde{u} - l^\top \tilde{v} + u^\top\tilde{w} + \bar{z} \sum_{i\colon \lambda_i(\tilde{S}) < 0} \lambda_i(\tilde{S}).
\end{equation*}
\end{theorem}
\ifJOC
\proof{Proof.}
\else
\begin{proof}
\fi
Let $Z^*$ be an optimal solution of~\eqref{eq:SDP} with the additional constraints $\Acal(Z) = b$ and $l \le \Bcal(Z) \le u$ and $(\tilde{y},\tilde{z},\tilde{u},\tilde{v},\tilde{w},\tilde{P})$ dual feasible. Then 
\begin{align*}
  \inprod{-W}{Z^*} &- ( \sum_{i=1}^n \tilde{y}_i + k\tilde{y}_{n+1} + b^\top \tilde{u} - l^\top \tilde{v} + u^\top\tilde{w}) \\
                   & = \inprod{-W}{Z^*} - \sum_{i=1}^n \tilde{y}_i\inprod{E_i}{Z^*} - \tilde{z}\inprod{I}{Z^*} - \inprod{\Acal(Z^*)}{\tilde{u}} + \inprod{\Bcal(Z^*)}{\tilde{v}} - \inprod{\Bcal(Z^*)}{\tilde{w}} \\
  &= \inprod{-W -  \sum_{i=1}^n \tilde{y}_iE_i - \tilde{y}_{n+1}I - \Acal^\top(\tilde{u}) + \Bcal^\top(\tilde{v}) - \Bcal^\top(\tilde{w})}{Z^*} \\
    & = \inprod{\tilde{P} + \tilde{S}}{Z^*} = \inprod{\tilde{P}}{Z^*}  + \inprod{\tilde{S}}{Z^*}.
\end{align*}
We have $\tilde{P}\ge 0$, $Z^* \ge 0$. Furthermore, the eigenvalues of $Z^*$ are nonnegative and bounded by one (Lemma~\ref{lem:eigboundZ}). Using this and Lemma~\ref{lem:jansson}, we obtain
\begin{align*}
  p^* = \inprod{-W}{Z^*} &\ge \sum_{i=1}^n \tilde{y}_i + k\tilde{y}_{n+1} + b^\top \tilde{u} - l^\top \tilde{v} + u^\top\tilde{w} + \inprod{\tilde{S}}{Z^*} \\
  & \ge \sum_{i=1}^n \tilde{y}_i + k\tilde{y}_{n+1} + b^\top \tilde{u} - l^\top \tilde{v} + u^\top\tilde{w} + \sum_{i\colon \lambda_i(\tilde{S}) < 0} \lambda_i(\tilde{S}).
\end{align*}
\ifJOC
\hfill\Halmos
\endproof
\else
\end{proof}
\fi


Before stating the result used in the branch-and-bound tree after merging data points, we introduce the following notation. Let $E_i^\ell$ be the symmetric matrix such that $\inprod{E_i^\ell}{Z^\ell} = (Z^\ell e^\ell)_i$.
\begin{corollary}
 Consider the SDP~\eqref{eq:SDPell} together with equations $\Acal(Z^\ell)=b$ (e.g., from cannot-link constraints) and inequalities $l \le \Bcal(Z^\ell) \le u$ (representing cutting planes) with optimal objective function value $p^*$. 
 Let $\tilde{y} \in \R^{n-\ell+1}$, $\tilde{u}$, $\tilde{v}$, $\tilde{w}$ being vectors of appropriate size, $\tilde{P}\in \snml$, $\tilde{P} \ge 0$ and set
 $\tilde{S} = -W^\ell - \sum_{i=1}^{n-\ell} \tilde{y}_i E^\ell_i + \tilde{y}_{n+\ell+1}\Diag(e^\ell)  - \Acal^\top(\tilde{u}) + \Bcal^\top(\tilde{v}) - \Bcal^\top(\tilde{w}) - \tilde{P}$. 
Then
\begin{equation*}
    p^* \ge \sum_{i=1}^{n-\ell} \tilde{y}_i + k\tilde{y}_{n-\ell+1}  + b^\top \tilde{u} - l^\top \tilde{v} + u^\top\tilde{w} + \sum_{i\colon \lambda_i(\tilde{S}) < 0} \lambda_i(\tilde{S}).
\end{equation*}

\end{corollary}
\ifJOC
\proof{Proof.}
\else
\begin{proof}
\fi
  Constraint~\eqref{eq:SDPellb} implies that the row-sum of any row in $Z^\ell$ is bounded by one since
  \begin{equation*}
      \sum_{j=1}^{n-\ell} z^\ell_{ij} \le \sum_{j=1}^{n-\ell} z^\ell_{ij}e^\ell_j = 1 \quad \mbox{for all } i \in \{1,\dots, n-\ell\}.
  \end{equation*}
  Hence using the same arguments as in Lemma~\ref{lem:eigboundZ} we can bound the eigenvalues by one and apply Theorem~\ref{thm:errorbound}.
\ifJOC
\hfill\Halmos
\endproof
\else
\end{proof}
\fi

\section{Heuristic}\label{sec:heuristic}

The most popular heuristic for solving MSSC is $k$-means \citep{macqueen1967some, lloyd1982least}. It can be viewed as a greedy algorithm. 
During each update step, all the data points are assigned to their nearest centers. Afterwards, the cluster centers are repositioned by calculating the mean of the assigned observations to the respective centroids. The update process is performed until the centroids are no longer updated and therefore all observations remain at the assigned clusters.
In this paper, we use COP $k$-means \citep{wagstaff2001constrained}, a constrained version of $k$-means that aims at finding high quality clusters using prior knowledge. COP $k$-means is a constrained clustering algorithm that belongs to a class of semi-supervised machine learning algorithms. Constrained clustering incorporates a set of must-link and cannot-link constraints that define a relationship between two data instances: a must-link constraint (ML) is used to specify that the two points in the must-link relation should be in the same cluster, whereas a cannot-link constraint (CL) is used to specify that the two points in the cannot-link relation should not be in the same cluster. These sets of constraints, which are naturally available as branching decisions while visiting the branch-and-bound tree, represent the prior knowledge on the problem for which $k$-means will attempt to find clusters that satisfy the specified ML and CL constraints. The algorithm returns an empty partition if no such clustering exists which satisfies the constraints. COP $k$-means is described in Algorithm~\ref{alg:copkmeans}. 

\begin{algorithm}
\caption{COP $k$-means}
\label{alg:copkmeans}

\FuncSty{K-MEANS(}\ArgSty{dataset $\mathcal{D}$, initial cluster centers $m_1, \dots, m_k$, must-link constraints ML $\subseteq \mathcal{D} \times \mathcal{D}$, cannot-link constraints CL $\subseteq \mathcal{D} \times \mathcal{D}$}\FuncSty{)}

\Repeat{convergence}{
\ForEach{data point $s_i \in \mathcal{D}$}{
$j \leftarrow \argmin \big\{ \|s_i - m_j\|^2 \colon j \in \{1,\dots,k\} \And $ \\ \hspace{3cm}$\texttt{VIOLATE\_CONSTRAINT}(s_i, C_j, \ML, \CL) \mathrm{~is~ false}\big\}$\;
\eIf{$j<\infty$}{
assign $s_i$ to $C_j$\;}{
\KwRet empty partition\;}} 
\ForEach{cluster $C_j$}{
$m_j \leftarrow $ mean of the data points $s_i$ assigned to $C_j$\;}
}
\KwRet $C_1, \dots, C_k$

\bigskip
\FuncSty{VIOLATE\_CONSTRAINTS(}\ArgSty{data point $s_i$, cluster $C_j$, must-link constraints $\ML \subseteq \mathcal{D} \times \mathcal{D}$, cannot-link constraints $\CL \subseteq \mathcal{D} \times \mathcal{D}$}\FuncSty{)}

\ForEach{$(s_i, s_h) \in \ML$}{
\lIf{$s_h \notin C_j$}{
\KwRet true}}
\ForEach{$(s_i, s_h) \in \CL$}{
\lIf{$s_h \in C_j$}{
\KwRet true}}
\KwRet false\;

\end{algorithm}

Like other local solvers for non-convex optimization
problems, $k$-means (both in the unconstrained and constrained version) is very sensitive to the choice of the initial centroids, therefore, it often converges to a local minimum rather than the global minimum of the MSSC objective. To overcome this drawback, the algorithm is initialized with several different starting points, choosing then the clustering with the lowest objective function \citep{franti2019much}.

In the literature, several initialization algorithms have been proposed to prevent $k$-means to get stuck in a low quality local minimum. The most popular strategy for initializing $k$-means is $k$-means++ \citep{arthur2006k}. The basic idea behind this approach is to spread out the $k$ initial cluster centers to avoid the poor clustering that can be found by the standard $k$-means algorithm with random initialization. More in detail, in $k$-means++, the first cluster center is randomly chosen from the data points. Then, each subsequent cluster center is chosen from the remaining data points with probability proportional to its squared distance from the already chosen cluster centers.

We aim to exploit the information available in the solution of the SDP relaxation in order to extract a centroid initialization for COP $k$-means.

In the literature, theoretical properties of the Peng-Wei relaxation have been studied under specific stochastic models. A feasible clustering can be derived by the solution of the SDP relaxation~(\ref{eq:SDP}) by a rounding step. Sometimes, the rounding step is unnecessary because the SDP relaxation finds a solution that is feasible for the original MSSC. This phenomenon is known in the literature as exact recovery or tightness of the relaxation. Recovery guarantees have been established under  a model called the subgaussian mixtures model, whose special cases include the stochastic ball model and Gaussian mixture model \citep{awasthi2015relax, iguchi2017probably, mixon2017clustering, li2020birds}. Under this distributional setting, cluster recovery is guaranteed with high probability whenever the distances between the clusters are sufficiently large. However, the generative assumption may not be satisfied by real data, and this implies that in general a rounding procedure is needed, and if possible also a bound improvement. Instead of building a rounding procedure, we decide to derive a ``smart'' initialization for the constrained $k$-means based on the solution of our bounding problem.
Here, we build the initialization exploiting the matrix $Z_{SDP}$ solution of the current bounding problem. The idea is that if the relaxation were tight, then $Z_{SDP}$ would be a clustering feasible for the rank constrained SDP \eqref{eq:RankSDP}, and hence would allow to easily recover the centroids. If the relaxation is not tight, the closest rank-$k$ approximation is built and it is used to recover the centroids. More in detail, let $Z$ be a feasible solution of the rank constrained SDP~\eqref{eq:RankSDP}. It is straightforward~\citep{mixon2017clustering} to see that $Z$ can be written as the sum of $k$ rank-one matrices:
\begin{equation}
\label{eq:clustering_matrix}
Z = \sum_{j=1}^{k} \frac{1}{|C_j|} \mathbbm{1}_{C_j} \mathbbm{1}^\top_{C_j},
\end{equation}
where $\mathbbm{1}_{C_j} \in \{0,1\}^n$ is the indicator vector of the $j$-th cluster, i.e., the $i$-th component of $\mathbbm{1}_{C_j}$ is 1 if the data point $p_i \in C_j$ and 0 otherwise. 
If we post-multiply $Z$ by the data matrix $W_p \in \mathbf{R}^{n\times d}$ whose $i$-th row is the data point $p_i$, we obtain a matrix  $M = ZW_p$ with a well defined structure. In fact, from equation (\ref{eq:clustering_matrix}) it follows that, for each $j \in \{1, \dots, k\}$, $M$ contains $|C_j|$ rows equal to the centroid of the data points assigned to $C_j$. If the SDP relaxation is tight, the different rows of $M$ are equal to the optimal centroids. In this case, it is natural to use the convex relaxation directly to obtain the underlying ground truth solution without the need for a rounding step. In practice, the optimizer of the SDP relaxation may not always be a clustering matrix, i.e., a low-rank solution as described by equation~(\ref{eq:clustering_matrix}).

The idea now is to build the rank-$k$ approximation $\hat{Z}$ which is obtained exploiting the following result.
\begin{proposition} \citep{eckart1936approximation}
Let $X$ be a positive semidefinite matrix with the eigenvalues $\lambda_1 \geq \lambda_2 \geq \ldots \geq \lambda_n\ge 0$ and the corresponding eigenvectors $v_1, v_2,\ldots,v_n$. If $X$ has rank $r$, for any $k < r$, the best rank $k$ approximation of $X$, for both the Frobenius and the spectral norms is given by 
\begin{equation}\label{eq:rankkapp}
\hat{X} = \sum_{i=i}^k \lambda_i v_i v_i^\top,
\end{equation}
which is the truncated eigenvalue decomposition of $X$.
\end{proposition}
Then, we compute the approximate centroid matrix $M=\hat{Z}W_p$.
In order to derive the $k$ centroids, the unconstrained $k$-means is applied to the rows of matrix $M$. Finally, the obtained centroids are used in order to initialize the algorithm COP $k$-means, which is run just once. The procedure is summarized in Algorithm~\ref{alg:sdpinit}.

\begin{algorithm}
\caption{SDP-based initialization of $k$-means}
\label{alg:sdpinit}

\FuncSty{SDP-INIT(}\ArgSty{dataset $\mathcal{D}$, number of clusters $k$, must-link constraints ML $\subseteq \mathcal{D} \times \mathcal{D}$, cannot-link constraints CL $\subseteq \mathcal{D} \times \mathcal{D}$}\FuncSty{)}

solve the SDP relaxation and obtain the optimizer $Z_{SDP}$\;
find the best rank $k$ approximation of $Z_{SDP}$ and obtain $\hat{Z}$ by \eqref{eq:rankkapp}\;
compute $M = \hat{Z}W_p$\;
cluster the rows of $M$ with unconstrained $k$-means to get the centroids $m_1, \dots, m_k$\;
use $m_1, \dots, m_k$ as the starting point of constrained $k$-means\;

\end{algorithm}

\begin{figure}

\begin{subfigure}{.5\linewidth}
 
\begin{tikzpicture}
\begin{axis}[
    title={LB = 1322.5, UB = 1384.5},
      width=\textwidth,
      height=\axisdefaultheight
]
\addplot[only marks, mark=o, mark size=2.5pt, color=gray]
table[]
{150_4_1.5.txt};

\addplot[only marks, mark=x, mark size=2.5pt]
table[]
{3/M1.txt};

\addplot[only marks, mark=diamond, color=black, mark size=4.0pt]
table[]
{3/C1.txt};
\end{axis}
\end{tikzpicture}
 
\end{subfigure}%
\begin{subfigure}{.5\linewidth}
 
\begin{tikzpicture}
\begin{axis}[
    title={LB = 1367.7, UB = 1384.5},
      width=\textwidth,
      height=\axisdefaultheight
]
\addplot[only marks, mark=o, mark size=2.5pt, color=gray]
table[]
{150_4_1.5.txt};

\addplot[only marks, mark=x, mark size=2.5pt]
table[]
{3/M2.txt};

\addplot[only marks, mark=diamond, color=black, mark size=4.0pt]
table[]
{3/C2.txt};
\end{axis}
\end{tikzpicture}
 
\end{subfigure}%

\begin{subfigure}{.5\linewidth}
 
\begin{tikzpicture}
\begin{axis}[
    title={LB = 1380.6, UB = 1383},
      width=\textwidth,
      height=\axisdefaultheight
]
\addplot[only marks, mark=o, mark size=2.5pt, color=gray]
table[]
{150_4_1.5.txt};

\addplot[only marks, mark=x, mark size=2.5pt]
table[]
{3/M3.txt};

\addplot[only marks, mark=diamond, color=black, mark size=4.0pt]
table[]
{3/C3.txt};
\end{axis}
\end{tikzpicture}
 
\end{subfigure}%
\begin{subfigure}{.5\linewidth}
 
\begin{tikzpicture}
\begin{axis}[
    title={LB = 1383, UB = 1383},
      width=\textwidth,
      height=\axisdefaultheight
]
\addplot[only marks, mark=o, mark size=2.5pt, color=gray]
table[]
{150_4_1.5.txt};

\addplot[only marks, mark=x, mark size=2.5pt]
table[]
{3/M4.txt};

\addplot[only marks, mark=diamond, color=black, mark size=4.0pt]
table[]
{3/C4.txt};
\end{axis}
\end{tikzpicture}
 
\end{subfigure}
\caption{An instance with 150~points and $k = 3$.}
\label{fig:heuristic3}
\end{figure}
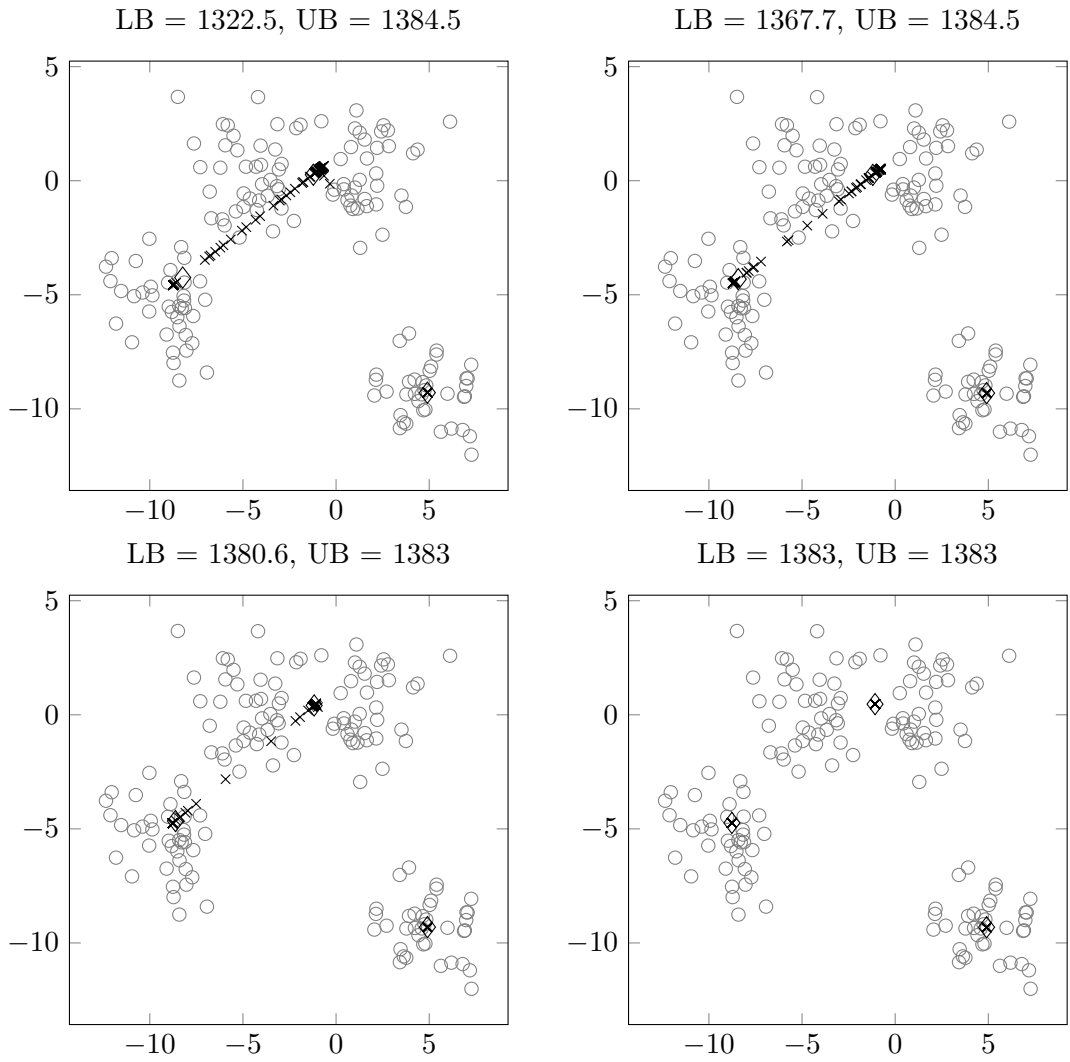

The intuition is that the better the SDP solution, the better the initialization, and hence the produced clustering. In order to confirm this intuition, we show the behavior of the heuristic on a synthetic example with 150 points in 2 dimensions. We denote by circles the points in $W_p$, by crosses the rows of matrix $M$ produced at Step 3, by diamonds the centroids obtained by clustering the rows of $M$ at Step 4 of the Algorithm~\ref{alg:sdpinit}. 
In Figure~\ref{fig:heuristic3} we assume $k=3$ and apply our heuristic on different solutions of the SDPs generated during our bounding procedure: in Figure~\ref{fig:heuristic3}~(a) we use as $Z_{SDP}$ the solution obtained by solving problem \eqref{eq:SDP}, and we can see that there is some gap (the upper and lower bounds are displayed on top of each figure) and that matrix $M$ has many different rows. In Figures~\ref{fig:heuristic3}~(b), (c), and (d) we consider as $Z_{SDP}$ the solution of the SDP obtained by performing respectively 1, 2 and 3 iterations of adding cutting-planes, i.e., solving problem~(\ref{eq:SDP}) with some additional constraints~(\ref{eq:pairs})--(\ref{eq:clique}). It is clear how the rows of $M$ converge to three different centroids that, in this case, correspond to the optimal solution (the gap here is zero).

The use of \texttt{SDP-INIT} as a standalone initialization procedure could be expensive since it needs to solve a certain number of SDP problems and to perform an eigenvalue decomposition on the solution that gives the best lower bound. However, when embedded in our branch-and-bound, the extra cost of running \texttt{SDP-INIT} is only the computation of the spectral decomposition of the SDP solution providing the lower bound at the node, which is negligible with respect to the bound computation. 

The effectiveness of the proposed heuristic algorithm is confirmed by the numerical results presented in Section~\ref{sec:heurnr}.

\section{Numerical Results}\label{sec:numericalresults}
In this section we describe the implementation details and we show the numerical results of \texttt{SOS-SDP} on synthetic and real-world datasets.

\subsection{Details on the Implementation}
\texttt{SOS-SDP} is implemented in C++ and we use as internal subroutine for computing the bound SDPNAL+ \citep{sdpnalplus, zhao-sun-toh-2010}, which is implemented in MATLAB. SDPNAL+ is called using the MATLAB Engine API that enables running MATLAB code from C++ programs. We note that solvers based on interior point methods are not practical when solving instances with such a large number of constraints. We run our experiments on a machine with Intel(R) Xeon(R) 8124M CPU @ 3.00GHz with 16 cores, 64 GB of RAM, and Ubuntu Server 20.04. The C++ Armadillo library \citep{sanderson2016armadillo} is extensively used to handle matrices and linear algebra operations efficiently.
\texttt{SOS-SDP} can be efficiently executed in a multi-thread environment. In order to guarantee an easy and highly configurable parallelization, we use the thread pool pattern. This pattern allows controlling the number of threads the branch-and-bound is creating and saving resources by reusing threads for processing different nodes of the tree.
We adopt the same branch-and-bound configuration for each instance. In particular, we visit the tree with the best-first search strategy. When the problem at a given level is divided into the \emph{must-link} and the \emph{cannot-link} sub-problems, each node is submitted to the thread pool and run in parallel with the other threads of the pool. 
Each thread of the branch-and-bound algorithm runs in a separate MATLAB session. Furthermore, since numerical algebra and linear functions are multi-threaded in MATLAB, these functions automatically execute on multiple computational threads in a single MATLAB session. To balance resource allocations for multiple MATLAB sessions and use all the available cores of the machine, we set a maximum number of computational threads allowed in each session.

\paragraph{Branch-and-bound setting}
On all the numerical tests, we adopt the following parameters setting. 
As for the pair and triangle inequalities, we randomly separate at most 100000 valid cuts, we sort them in decreasing order with respect to the violation, and we select the first 5\% of violated ones, yielding at most 5000 pairs and at most 5000 triangles added in each cutting-plane iteration. Since effective inequalities are inherited from the parent to its children, at the root node the maximum number of cutting-plane iterations is set to $cp_{\mathrm{max}} = 50$, whereas for the children this number is set to 30. 
The tolerance for checking the violation of the cuts is set to $\varepsilon_{\mathrm{viol}} = 10^{-4}$, whereas the tolerance for identifying the active inequalities is set to $\varepsilon_{\mathrm{act}} = 10^{-6}$. Finally, we set the accuracy tolerance of SDPNAL+ to $10^{-5}$. 

As for the parallel setting, we use different configurations depending on the size of the instances since the solver requires a higher number of threads to efficiently solve large size problems.
For small instances ($n < 500$) we create a pool of 16 threads, each of them running on a session with a single component thread. For medium instances ($500 \leq n < 1000$) we use a pool of 8 threads, each of them running on a session with 2 component threads. For ($1000 \leq n < 1500$) we use a pool of 4 threads, each of them runs on a session with 4 component threads. Finally, for large scale instances ($n \geq 1500$) we use a pool of 2 threads, each of them running on a session with 8 component threads. In all cases, the MATLAB session for the computation at the root node uses all the available cores. The source code is available at \url{https://github.com/INFORMSJoC/2021.0096} \citep{SOS-SDP2021}.

\subsection{Benchmark Instances}
In order to test extensively the efficiency of \texttt{SOS-SDP} we use both artificial datasets that are built in such a way to be compliant with the MSSC assumptions and real-world datasets.

\paragraph{Artificial Instances}
Due to the minimization of the sum of squared Euclidean distances, an algorithm that solves the MSSC finds spherically distributed clusters around the centers. In order to show the effectiveness of our algorithm on instances compliant with the MSSC assumptions, we generate very large scale Gaussian datasets in the plane $(d=2)$ with varying number of data points $n \in \{2000, 2500, 3000 \}$, number of clusters $k \in \{10, 15\}$ and degree of overlap. More in detail, we sample $n$ points from a mixture of $k$ Gaussian distributions $\mathcal{N}(\mu_j, \Sigma_j)$ with equal mixing proportions, mean $\mu_j$ and shared spherical covariance matrix $\Sigma_j = \sigma I$, where $\sigma \in \{0.5, 1.0\}$ is the standard deviation. The cluster centers $\mu_j$ are sampled from a uniform distribution in the interval $[-\frac{n}{1000}-k, \frac{n}{1000}+k]$. We use the following notation to name the instances: $\{n\}\_\{k\}\_\{\sigma\}$. Note that in this case, we know in advance the correct number of clusters, so we only solve the instances for that value of~$k$.

\paragraph{Real-world Datasets}
We use a set of 34 real-world datasets coming from different domains, with a number of entities $n$ ranging between $75$ and $4177$, and with a number of features $d$ ranging between $2$ and $20531$. The datasets' characteristics are reported in Table~\ref{tab:datasets}. 

\begin{table}
\begin{center}
\begin{tabular}{lcc}
\toprule
Dataset   &  $n$   &  $d$   \\
\midrule
Ruspini & 75 & 2 \\
Voice & 126 & 310 \\
Iris & 150 & 4 \\
Wine & 178 & 13 \\
Gr202 & 202 & 2 \\
Seeds & 210 & 7 \\
Glass & 214 & 9 \\
CatsDogs & 328 & 14773 \\
Accent & 329 & 12 \\
Ecoli & 336 & 7 \\
RealEstate & 414 & 5 \\
Wholesale & 440 & 11 \\
ECG5000 & 500 & 140 \\
Hungarian & 522 & 20 \\
Wdbc & 569 & 30 \\
Control & 600 & 60 \\
Heartbeat & 606 & 3053 \\
\bottomrule
\end{tabular}\qquad 
\begin{tabular}{lcc}
\toprule
Dataset   &  $n$   &  $d$   \\
\midrule
Strawberry & 613 & 235 \\
Energy & 768 & 16 \\
Gene & 801 & 20531 \\
SalesWeekly &  810 & 106 \\
Vehicle & 846 & 18 \\
Arcene & 900 & 10000 \\
Wafer & 1000 & 152 \\
Power & 1096 & 24 \\
Phishing & 1353 & 9 \\
Aspirin & 1500 & 63 \\
Car & 1727 & 11 \\
Wifi & 2000 & 7 \\
Ethanol & 2000 & 27 \\
Mallat & 2400 & 1024 \\
Advertising & 3279 & 1558 \\
Rice & 3810 & 7 \\
Abalone & 4177 & 10 \\
\bottomrule
\end{tabular}
\end{center}
\caption{Characteristics of the real world datasets. They all can be downloaded at the UCI \citep{uci}, UCR \citep{UCRArchive2018} and sGDML \citep{chmiela2019sgdml} websites.}
\label{tab:datasets}
\end{table}

\subsection{Branch-and-Bound Results on Artificial instances}

In Table \ref{tab:res_art} we report the dataset name according to the notation $\{n\}\_\{k\}\_\{\sigma\}$, the optimal objective function $f_\textrm{opt}$, the number of cutting-plane iterations at the root (cp), the number of cuts added in the last cutting-plane iteration at the root ($\textrm{cuts}_\textrm{cp}$), the gap at the root ($gap_0$) when problem \eqref{eq:SDPbab} is solved without adding valid inequalities, in brackets the gap at the end of the cutting-plane procedure at the root node ($gap_{cp}$),
the number of nodes of the branch-and-bound tree (N), and the wall clock time in seconds (time).

\begin{tabularx}{\textwidth}{lcccccccc}
\toprule
Dataset 	&	$f_\textrm{opt}$	&	cp	&	$\textrm{cuts}_\textrm{cp}$	&	$gap_0$ $(gap_{cp})$	&	N	&	time	\\
\midrule
   2000\_10\_0.5 &  955.668 &             0 &          0 & 0.000039 (0.000039) &        1 &    848.88 \\
   2000\_10\_1.0 & 3601.310 &             3 &      10999 & 0.006171 (0.003578) &        3 &   8794.17 \\
   2000\_15\_0.5 &  955.800 &             1 &       6177 & 0.001556 (0.000009) &        1 &   1155.06 \\
   2000\_15\_1.0 & 3658.730 &             3 &      11035 & 0.006192 (0.002059) &        3 &   8351.91 \\
   2500\_10\_0.5 & 1199.080 &             1 &       5249 & 0.000184 (0.000083) &        1 &   2859.30 \\
   2500\_10\_1.0 & 4522.350 &            12 &      11539 & 0.008008 (0.000553) &        1 &  20495.43 \\
   2500\_15\_0.5 & 1194.550 &             0 &          0 & 0.000699 (0.000699) &        1 &   1049.76 \\
   2500\_15\_1.0 & 4574.360 &             6 &      10146 & 0.005311 (0.000971) &        1 &  10245.69 \\
   3000\_10\_0.5 & 1446.480 &             0 &          0 & 0.000067 (0.000067) &        1 &   2220.21 \\
   3000\_10\_1.0 & 5512.370 &             9 &      10769 & 0.004601 (0.000606) &        1 &  27781.38 \\
   3000\_15\_0.5 & 1439.940 &             0 &          0 & 0.000433 (0.000433) &        1 &   2003.94 \\
   3000\_15\_1.0 & 5537.200 &            10 &      15608 & 0.006245 (0.001205) &        3 &  38330.01 \\
\bottomrule
\caption{Results for the artificial datasets.}
\label{tab:res_art}
\end{tabularx}

\begin{figure}
\begin{center} 
\begin{subfigure}{.5\linewidth}
 
\begin{tikzpicture}
\begin{axis}[
    title={$k=10, \ \sigma=0.5$},
      width=\textwidth,
      height=\axisdefaultheight
]
\addplot[only marks, mark=o, mark size=1.0pt]
table[]
{2000_10_0.5.txt};
\end{axis}
\end{tikzpicture}
 
\end{subfigure}%
\begin{subfigure}{.5\linewidth}
 
\begin{tikzpicture}
\begin{axis}[
    title={$k=10, \ \sigma=1.0$},
      width=\textwidth,
      height=\axisdefaultheight
]
\addplot[only marks, mark=o, mark size=1.0pt]
table[]
{2000_10_1.txt};
\end{axis}
\end{tikzpicture}
 
\end{subfigure}%

\begin{subfigure}{.5\linewidth}
 
\begin{tikzpicture}
\begin{axis}[
    title={$k=15, \ \sigma=0.5$},
      width=\textwidth,
      height=\axisdefaultheight
]
\addplot[only marks, mark=o, mark size=1.0pt]
table[]
{2000_15_0.5.txt};
\end{axis}
\end{tikzpicture}
 
\end{subfigure}%
\begin{subfigure}{.5\linewidth}
 
\begin{tikzpicture}
\begin{axis}[
    title={$k=15, \ \sigma=1.0$},
      width=\textwidth,
      height=\axisdefaultheight
]
\addplot[only marks, mark=o, mark size=1.0pt]
table[]
{2000_15_1.txt};
\end{axis}
\end{tikzpicture}
 
\end{subfigure}
\caption{Artificial instances for $n=2000$ and $d=2$.}
\label{fig:art}
\end{center}
\end{figure}
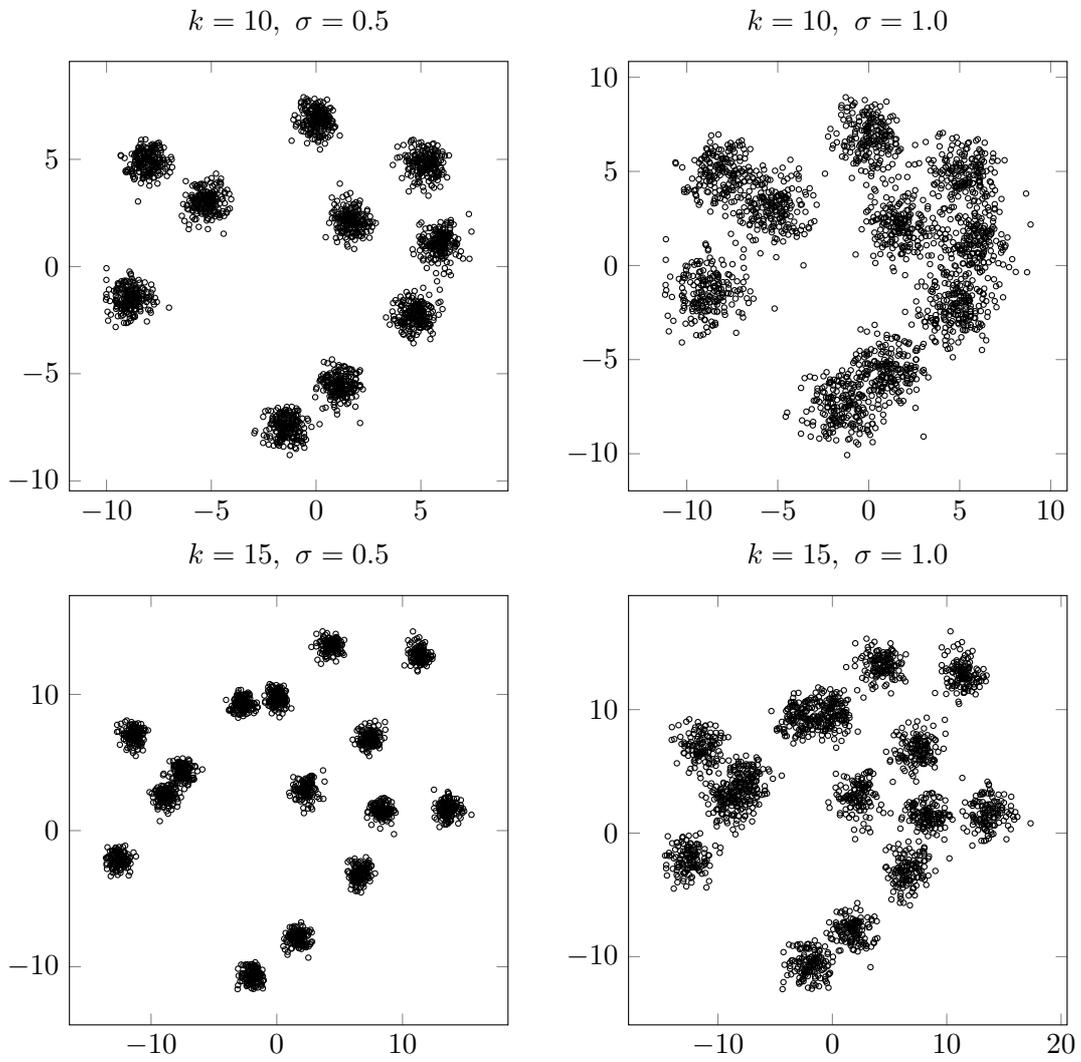

As we increase $\sigma$, the cluster separation decreases, and the degree of overlap increases (see Figure \ref{fig:art}). In this scenario, the SDP relaxation is not tight anymore and the global minimum is certified by our specialized branch-and-bound algorithm.
For $\sigma = 0.5$ each problem is solved at the root with zero (i.e., the SDP relaxation is tight) or with at most one cutting-plane iteration. As we decrease the cluster separation by increasing $\sigma$ the problem becomes harder since some clusters overlap and the cluster boundaries are less clear. In this case, more cutting-plane iterations are needed (up to a maximum of 12 iterations). In any case, we need at most 3 nodes for solving these instances, and this confirms that, if the generative assumption is met, the cutting-plane procedure at the root node is the main ingredient for success. In the next section, we show how the behavior changes in real world instances, where we do not have information on the data distribution and on the correct value of $k$. In this case, the overall branch-and-bound algorithm becomes fundamental in order to solve the problems.

\subsection{Branch-and-Bound Results on Real World Datasets}

The MSSC requires the user to specify the number of clusters $k$ to generate. Determining the right $k$ for a data set is a different issue from the process of solving the clustering problem. This is still an open problem since, depending on the chosen distance measure, one value of $k$ may be better than another one. Hence, choosing $k$ is often based on assumptions on the application, prior knowledge of the properties of the dataset, and practical experience. In the literature, clustering validity indices in conjunction with the $k$-means algorithm are commonly used to determine the ``right'' number of clusters. Most of these methods minimize or maximize an external validity index by running a clustering algorithm (for example $k$-means) several times for different values of $k$.
We recall that the basic idea behind the MSSC is to define clusters such that the total within-cluster sum of squares is minimized. This objective function measures the compactness of the clustering and we want it to be as small as possible. The ``elbow method'' is probably the most popular method for determining the number of clusters. It requires running the $k$-means algorithm with an increasing number of clusters. The suggested $k$ can be determined by looking at the MSSC objective as a function of $k$ and by finding the inflection point. The location of the inflection point (knee) in the plot is generally considered as an indicator of the appropriate number of clusters. The drawback of this method is that the identification of the knee could not be obvious. Hence, different validity indices have been proposed in the literature to identify the suitable number of clusters or to check whether a given dataset exhibits some kind of a structure that can be captured by a clustering algorithm for a given $k$. All these indices are computed aposteriori given the clustering produced for different values of $k$. In addition to the elbow method, we use three cluster validity measures that are compliant with the assumptions of the MSSC: namely the Silhouette index \citep{rousseeuw1987silhouettes}, the Calinski–Harabasz (CH) index \citep{calinski1974dendrite} and Davies–Bouldin (DB) index \citep{davies1979cluster}. 
The Silhouette index determines how well each object lies within its cluster and is given by the average Silhouette coefficient over all the data points. The Silhouette coefficient is defined for each data point and is composed of two scores: the mean distance between a sample and all other points in the same class and the mean distance between a sample and all other points in the next nearest cluster. The CH index is the ratio of the sum of between-clusters dispersion and within-cluster dispersion for all clusters. The DB index is defined as the average similarity between each cluster and its most similar one.

The Silhouette index and the CH index are higher when the clusters are dense and well separated, which relates to the standard concept of clustering, whereas for the DB index lower values indicate a better partition. 

Since the exact resolution of the MSCC problem could be expensive and time consuming from the computational point of view, one may be interested in finding the global solution for a specified or restricted number of clusters. In practice, one can run the $k$-means algorithm for different values of $k$ and then use the exact algorithm to find and certify the global optimum for the $k$ suitable for the application of interest.
Hence, we choose to run our algorithm on a large number of datasets, and for each dataset, we run it only for the suggested number of clusters obtained with the help of the criteria mentioned above. 
Whenever there is some ambiguity, i.e., the different criteria suggest different values of $k$, we run our algorithm for all the suggested values. With this criterion, we end up solving $54$ clustering instances with different size $n$, different dimension $d$, and different values of $k$.

In Table \ref{tab:res_dataset} we report:
\begin{itemize}
    \item the dataset name
    \item the number of clusters ($k$)
    \item the optimal objective function ($f_\textrm{opt}$). We add a $(*)$ whenever the optimum we certify is not found by $k$-means at the root node
    \item the number of cutting-plane iterations at root (cp)
    \item the number of inequalities of the last SDP problem solved at the root in the cutting-plane procedure ($\textrm{cuts}_\textrm{cp}$)
    \item the gap at the root ($gap_0$) when problem \eqref{eq:SDPbab} is solved without adding valid inequalities, and in brackets the gap at the end of the cutting-plane procedure at the root node ($gap_{cp}$) 
    \item the number of nodes (N) of the branch-and-bound tree
    \item the wall clock time in seconds (time).
\end{itemize}

Small and medium scale instances ($n < 1000$) are considered solved when the relative gap tolerance
is less or equal than $10^{-4}$, whereas for large scale instances ($n \geq 1000$) the branch-and-bound algorithm is stopped when the tolerance is less or equal than $10^{-3}$, which we feel is an adequate tolerance for large scale real-world applications. The gap measures the difference between the best upper and lower bounds and it is calculated as $(UB - LB) / UB$.

The numerical results show that our method is able to solve successfully all the instances up to a size of $n=4177$ entities. 


\begin{longtable}{lcccccccc}
\toprule
Dataset	&	$k$	&	$f_\textrm{opt}$	&	cp	&	$\textrm{cuts}_\textrm{cp}$	&	$gap_0$ $(gap_{cp})$	&	N	&	time	\\
\midrule
Ruspini	&	4	&	1.28811e+04	&	0	&	0	&	2.23e-04	 (2.23e-04)	&	1	&	2.55	\\
Voice	&	2	&	1.13277e+22	&	2	&	7593	&	5.40e-02	 (1.66e-06)	&	1	&	14.45	\\
Voice	&	9	&	5.74324e+20*	&	4	&	6115	&	1.07e-01	 (6.45e-04)	&	3	&	128.35	\\
Iris	&	2	&	1.52348e+02	&	2	&	7701	&	1.10e-02	 (2.19e-06)	&	1	&	17	\\
Iris	&	3	&	7.88514e+01	&	4	&	7136	&	4.23e-02	 (1.18e-04)	&	5	&	83.3	\\
Iris	&	4	&	5.72285e+01	&	4	&	7262	&	4.28e-02	 (4.20e-04)	&	3	&	104.55	\\
Wine	&	2	&	4.54375e+06	&	3	&	8162	&	3.45e-02	 (2.69e-07)	&	1	&	53.55	\\
Wine	&	7	&	4.12138e+05*	&	4	&	5759	&	5.81e-02	 (1.03e-04)	&	3	&	87.55	\\
Gr202	&	6	&	6.76488e+03	&	6	&	6607	&	6.72e-02	 (8.53e-04)	&	17	&	298.35	\\
Seeds	&	2	&	1.01161e+03*	&	9	&	10186	&	4.31e-02	 (4.77e-04)	&	29	&	957.1	\\
Seeds	&	3	&	5.87319e+02	&	4	&	6620	&	2.67e-02	 (1.26e-05)	&	1	&	68.85	\\
Glass	&	3	&	1.14341e+02	&	5	&	6799	&	4.68e-02	 (1.64e-04)	&	3	&	193.8	\\
Glass	&	6	&	7.29647e+01*	&	7	&	3014	&	5.45e-02	 (4.36e-04)	&	5	&	198.9	\\
CatsDogs	&	2	&	1.14099e+05	&	1	&	5368	&	1.83e-03	 (2.23e-09)	&	1	&	108.8	\\
Accent	&	2	&	3.28685e+04	&	0	&	0	&	6.55e-06	 (6.55e-06)	&	1	&	11.05	\\
Accent	&	6	&	1.84360e+04*	&	8	&	4523	&	2.94e-02	 (2.08e-05)	&	1	&	244.8	\\
Ecoli	&	3	&	2.32610e+01	&	4	&	10101	&	7.71e-03	 (1.89e-04)	&	3	&	181.9	\\
RealEstate	&	3	&	5.50785e+07	&	3	&	6236	&	1.59e-02	 (3.51e-05)	&	1	&	104.55	\\
RealEstate	&	5	&	2.18711e+07	&	5	&	8006	&	6.82e-02	 (2.64e-05)	&	1	&	258.4	\\
Wholesale	&	5	&	2.04735e+03	&	6	&	7668	&	6.43e-02	 (2.06e-05)	&	1	&	421.6	\\
Wholesale	&	6	&	1.73496e+03*	&	10	&	11161	&	6.32e-02	 (7.06e-04)	&	3	&	1782.45	\\
ECG5000	&	2	&	1.61359e+04	&	3	&	9312	&	1.02e-03	 (7.49e-05)	&	1	&	119	\\
ECG5000	&	5	&	1.15458e+04	&	25	&	6289	&	4.93e-02	 (1.01e-04)	&	3	&	2524.5	\\
Hungarian	&	2	&	8.80283e+06	&	7	&	11265	&	1.03e-02	 (1.31e-05)	&	1	&	551.65	\\
Wdbc	&	2	&	7.79431e+07	&	5	&	8645	&	3.21e-02	 (2.10e-05)	&	1	&	436.05	\\
Wdbc	&	5	&	2.05352e+07*	&	23	&	10662	&	7.45e-02	 (5.27e-04)	&	15	&	2436.95	\\
Control	&	3	&	1.23438e+06	&	6	&	12381	&	2.80e-03	 (1.26e-04)	&	9	&	895.9	\\
Heartbeat	&	2	&	2.79391e+04	&	0	&	0	&	8.15e-06	 (8.15e-06)	&	1	&	66.3	\\
Strawberry	&	2	&	2.79363e+03	&	15	&	23776	&	5.44e-02	 (4.02e-04)	&	37	&	5250.45	\\
Energy	&	2	&	9.64123e+03	&	0	&	0	&	9.86e-09	 (9.86e-09)	&	1	&	18.7	\\
Energy	&	12	&	4.87456e+03	&	0	&	0	&	4.03e-07	 (4.03e-07)	&	1	&	29.75	\\
Gene	&	5	&	1.78019e+07*	&	2	&	15589	&	1.83e-03	 (1.30e-04)	&	3	&	3851.35	\\
Gene	&	6	&	1.70738e+07	&	5	&	14620	&	3.82e-03	 (2.08e-04)	&	11	&	9896.55	\\
SalesWeekly	&	2	&	1.44942e+06*	&	6	&	8508	&	2.50e-02	 (1.33e-03)	&	9	&	2341.75	\\
SalesWeekly	&	3	&	7.09183e+05*	&	4	&	9096	&	1.03e-03	 (9.44e-05)	&	1	&	262.65	\\
SalesWeekly	&	5	&	5.20938e+05*	&	4	&	11811	&	1.67e-03	 (1.12e-04)	&	5	&	1045.5	\\
Vehicle	&	2	&	7.29088e+06	&	5	&	10395	&	7.68e-03	 (3.72e-04)	&	11	&	1842.8	\\
Arcene	&	2	&	3.48490e+10	&	3	&	36100	&	2.59e-03	 (1.26e-04)	&	3	&	1369.35	\\
Arcene	&	3	&	2.02369e+10	&	0	&	0	&	3.50e-06	 (3.50e-06)	&	1	&	758.2	\\
Arcene	&	5	&	1.69096e+10*	&	7	&	8327	&	7.57e-03	 (1.55e-04)	&	27	&	6885	\\
Wafer	&	2	&	6.19539e+04	&	3	&	7254	&	7.82e-04	 (1.00e-04)	&	1	&	379.1	\\
Wafer	&	4	&	4.42751e+04	&	22	&	16957	&	1.97e-02	 (8.76e-04)	&	1	&	6756.65	\\
Power &	2	&	3.22063e+03	&	3	&	11350	&	1.05e-02	 (2.89e-03)	&	3	&	3381.3	\\
Phishing	&	9	&	3.15888e+03*	&	46	&	12459	&	2.48e-02	 (7.00e-04)	&	1	&	18866.6	\\
Aspirin	&	3	&	1.27669e+04	&	2	&	10000	&	4.39e-03	 (3.02e-03)	&	9	&	3779.1	\\
Car	&	4	&	5.61600e+03	&	23	&	38582	&	1.61e-03	 (1.02e-05)	&	1	&	5989.95	\\
Ethanol	&	2	&	7.26854e+03	&	0	&	0	&	5.33e-08	 (5.33e-08)	&	1	&	310.25	\\
Wifi	&	5	&	2.04311e+05	&	7	&	20886	&	1.13e-02	 (2.18e-03)	&	7	&	22754.5	\\
Mallat	&	3	&	9.08648e+04	&	5	&	17092	&	3.61e-03	 (9.59e-04)	&	1	&	5970.4	\\
Mallat	&	4	&	7.45227e+04	&	6	&	15305	&	6.80e-03	 (4.49e-03)	&	5	&	26344.9	\\
Advertising	&	2	&	5.00383e+06*	&	1	&	12533	&	1.53e-03	 (2.16e-05)	&	1	&	6465.1	\\
Advertising	&	8	&	4.54497e+06*	&	4	&	19948	&	2.98e-03	 (1.08e-04)	&	1	&	25114.1	\\
Rice	&	2	&	1.39251e+04	&	24	&	7258	&	1.43e-02	 (7.14e-03)	&	5	&	103710.2	\\
Abalone	&	3	&	1.00507e+03	&	0	&	0	&	3.14e-04	 (3.14e-04)	&	1	&	9428.2	\\
\bottomrule
\caption{Results for the real world datasets}
\label{tab:res_dataset}
\end{longtable}


To the best of our knowledge, the exact algorithm proposed in \cite{aloise2012improved} represents the actual state-of-the-art. Indeed it is the only algorithm able to exactly solve instances of size larger than 1000, satisfying one of the following strong assumptions (due to the geometrical approach involved): either the instance is on the plane ($d=2$) or the required number of clusters is large with respect to the number of points. Indeed they were able to solve a TSP instance with $d=2$ of size $n=2392$ for numbers of clusters ranging from $k=2$ to $k=10$, and for large number of clusters ($k$ between 100 and 400), and an instance of size $n=2310$ with $d=19$ but only for large number of clusters ($k$ between 230 and 500). Our algorithm has orthogonal capabilities in some sense to the one proposed in \cite{aloise2012improved}, since is not influenced by the number of features (we solve problems with thousands of features, which would be completely out of reach for the algorithm in~\cite{aloise2012improved}). Indeed, in the SDP formulation, the number of features is hidden in the matrix $W$, which is computed only once, so that it does not influence the computational cost of the algorithm. On the other hand, it is well known that the difficulty (and the gap) of the SDP relaxation \eqref{eq:SDPbab} increases when the boundary of the clusters are confused, and this phenomenon becomes more frequent when the number of clusters is high with respect to the number of points, and far away from the correct $k$ for the MSSC objective function. The strength of our bounding procedure is confirmed by 28 problems out of 54 solved at the root. Among these 28 problems, only 8 are tight, in the sense that problem \eqref{eq:SDPbab} without inequalities produces the optimal solution.
The efficiency of \texttt{SOS-SDP} comes from the combination of the cutting-plane procedure that allows us to close a significant amount of the gap even when the bound without inequalities is not tight, and the heuristic that when the SDP solution is good allows us to find the optimal solution. Note that in 15 out of 34 instances, our algorithm certifies the optimality of a solution that $k$-means at the root could not find.

Overall, the number of nodes of the branch-and-bound tree is always smaller than 40, but the computational cost of the single node may be high due to the high number of cutting-plane iterations. The values of cuts$_{\mathrm{cp}}$ confirm that the removal of inactive inequalities is effective, and allows to keep the number of inequalities moderate so that the SDP at each cutting-plane iteration is computationally tractable.

\subsection{Numerical Results of \texttt{SDP-INIT} }\label{sec:heurnr}

In order to test the efficiency of our initialization of constrained $k$-means, we report the behaviour at the root node on a subset of real-world datasets. We selected the most popular on the UCI website with size in the range of 150--569.  To have more difficult instances, we run the heuristic for all the values of $k$ in the range from $2$ to $10$. Note that for $k$ far from the values suggested by the validation indices, the optimal solution may be constituted by overlapped and confused clusters that are more difficult to find for any heuristic. 

In Table~\ref{tab:heuristic}, we report the results obtained by our heuristic, compared with 50 runs of $k$-means initialized with $k$-means++ and with random initialization.

In each table, we report:
\begin{itemize}
    \item the lower bound obtained by solving the basic SDP relaxation ($LB_0$), and the corresponding heuristic solution ($UB_0$)
    \item the lower bound obtained after performing $CP$ cutting-plane iterations $LB_{CP}$ and the corresponding heuristic solution ($UB_{CP}$)
    \item the solution produced by $k$-means after 50 runs initialized with $k$-means++ ($UB_{++}$)
    \item the solution produced by $k$-means after 50 runs randomly initialized ($UB_{RAND}$)
\end{itemize} 
We highlight the best solution in boldface. The results show that the solution $UB_{CP}$ is always the best, apart from 1 case. Note that in many cases, the solution $UB_{0}$ is fairly competitive both in terms of bound quality and computational effort since it requires the solution of exactly one SDP.

\begin{table}
\begin{center}
\footnotesize
\begin{tabular}{cccccccc}
\toprule
 $K$ &  $CP$ &     $LB_{0}$ &       $LB_{CP}$ &  $UB_{0}$ &   $UB_{CP}$ &      $UB_{++}$ &    $UB_{RAND}$ \\
\midrule
\multicolumn{8}{l}{Iris dataset}\\
 \midrule
 2 &        2 & 1.50679e+02 & 1.52348e+02 &   \bftab{1.52348e+02} & \bftab{1.52348e+02} & \bftab{1.52348e+02} & \bftab{1.52348e+02} \\
 3 &        4 & 7.55144e+01 & 7.88421e+01 &   7.88557e+01 & \bftab{7.88514e+01} & 7.88518e+01 & 7.88527e+01 \\
 4 &        6 & 5.47766e+01 & 5.72281e+01 &   \bftab{5.72285e+01} & \bftab{5.72285e+01} & 5.72560e+01 & 5.72560e+01 \\
 5 &        3 & 4.38467e+01 & 4.64369e+01 &   4.64612e+01 & \bftab{4.64462e+01} & \bftab{4.64462e+01} & 4.64612e+01 \\
 6 &        4 & 3.67110e+01 & 3.90175e+01 &   3.90660e+01 & \bftab{3.90400e+01} & 3.90660e+01 & \bftab{3.90400e+00} \\
 7 &        6 & 3.18467e+01 & 3.42788e+01 &   3.43058e+01 & \bftab{3.42982e+01} & 3.44090e+01 & 3.43859e+01 \\
 8 &        3 & 2.88697e+01 & 2.99660e+01 &   \bftab{2.99904e+01} & \bftab{2.99904e+01} & \bftab{2.99904e+01} & 3.04762e+01 \\
 9 &        3 & 2.64849e+01 & 2.77836e+01 &   2.79408e+01 & \bftab{2.77861e+01} & 2.78921e+01 & 2.83071e+01 \\
10 &        4 & 2.44186e+01 & 2.58329e+01 &   2.62712e+01 & \bftab{2.58341e+01} & 2.59644e+01 & 2.65776e+01 \\
\midrule 
\multicolumn{8}{l}{Glass dataset}\\
\midrule 
 2 &        6 & 1.35499e+02 & 1.36525e+02 &   1.36537e+02 & \bftab{1.36528e+02} & \bftab{1.36528e+02} & 1.36537e+02 \\
 3 &        5 & 1.08991e+02 & 1.14320e+02 &   \bftab{1.14341e+02} & \bftab{1.14341e+02} & \bftab{1.14341e+02} & \bftab{1.14341e+02} \\
 4 &        5 & 9.14749e+01 & 9.47742e+01 &   9.48402e+01 & \bftab{9.47899e+01} & 9.48402e+01 & 9.48402e+01 \\
 5 &        6 & 7.87104e+01 & 8.34045e+01 &   8.40062e+01 & \bftab{8.35054e+01} & 8.42973e+01 & 8.40502e+01 \\
 6 &        8 & 6.89918e+01 & 7.29430e+01 &   \bftab{7.29647e+01} & \bftab{7.29647e+01} & 7.37947e+01 & 7.43696e+01 \\
 7 &        8 & 6.19552e+01 & 6.47908e+01 &   6.53398e+01 & \bftab{6.47973e+01} & 7.08087e+01 & 6.66828e+01 \\
 8 &        6 & 5.61534e+01 & 5.85654e+01 &   5.87606e+01 & \bftab{5.85699e+01} & 5.90119e+01 & 6.08941e+01 \\
 9 &       10 & 5.12932e+01 & 5.37277e+01 &   5.41810e+01 & \bftab{5.37580e+01} & 5.55979e+01 & 5.61847e+01 \\
10 &        4 & 4.70718e+01 & 4.93411e+01 &   4.97866e+01 & \bftab{4.97382e+01} & 5.15837e+01 & 5.25047e+01 \\
\midrule
\multicolumn{8}{l}{Wholesale dataset}\\
\midrule 
 2 &        2 & 3.48221e+03 & 3.48656e+03 &   \bftab{3.48657e+03} & \bftab{3.48657e+03} & \bftab{3.48657e+03} & \bftab{3.48657e+03} \\
 3 &        5 & 2.85705e+03 & 2.91234e+03 &   \bftab{2.91252e+03} & \bftab{2.91252e+03} &\bftab{2.91252e+03} & \bftab{2.91254e+03} \\
 4 &        9 & 2.33207e+03 & 2.46555e+03 &   \bftab{2.46558e+03} & \bftab{2.46558e+03} & \bftab{2.46558e+03} & \bftab{2.46558e+03} \\
 5 &        7 & 1.91575e+03 & 2.04735e+03 &   2.04741e+03 & \bftab{2.04735e+03} & 2.04891e+03 & 2.04891e+03 \\
 6 &       10 & 1.63098e+03 & 1.73382e+03 &   1.74322e+03 & \bftab{1.73496e+03} & 1.74096e+03 & 1.75359e+03 \\
 7 &       12 & 1.44236e+03 & 1.52350e+03 &   1.52551e+03 & \bftab{1.52383e+03} & 1.52693e+03 & 1.53677e+03 \\
 8 &       11 & 1.28695e+03 & 1.36289e+03 &   1.36949e+03 & \bftab{1.36290e+03} & 1.36621e+03 & 1.39735e+03 \\
 9 &       10 & 1.14692e+03 & 1.21928e+03 &   1.22008e+03 & \bftab{1.21978e+03} & \bftab{1.21978e+03} & 1.26105e+03 \\
10 &        6 & 1.03078e+03 & 1.07843e+03 &   1.08010e+03 & \bftab{1.07843e+03} & 1.13670e+03 & 1.21282e+03 \\
\midrule
\multicolumn{8}{l}{Wdbc dataset}\\
\midrule 
 2 &        6 & 7.54429e+07 & 7.79415e+07 &   \bftab{7.79431e+07} & \bftab{7.79431e+07} & \bftab{7.79431e+07} & \bftab{7.79431e+07} \\
 3 &       27 & 4.14673e+07 & 4.72612e+07 &   4.74219e+07 & \bftab{4.72648e+07} & \bftab{4.72648e+07} & 4.74999e+07 \\
 4 &       22 & 2.62662e+07 & 2.91013e+07 &   2.92269e+07 & \bftab{2.92265e+07} & \bftab{2.92265e+07} & \bftab{2.92265e+07} \\
 5 &       20 & 1.90062e+07 & 2.05248e+07 &   2.05806e+07 & \bftab{2.05352e+07} & \bftab{2.05352e+07} & 2.06727e+07 \\
 6 &        6 & 1.47880e+07 & 1.55897e+07 &   1.69771e+07 & 1.69343e+07 & \bftab{1.66461e+07} & 1.71215e+07 \\
 7 &       22 & 1.20747e+07 & 1.31868e+07 &   1.32742e+07 & \bftab{1.32470e+07} & 1.32655e+07 & 1.33533e+07 \\
 8 &        8 & 1.02027e+07 & 1.07390e+07 &   1.12114e+07 & \bftab{1.12064e+07} & 1.12441e+07 & 1.15090e+07 \\
 9 &        3 & 8.83658e+06 & 9.09983e+06 &   \bftab{9.43290e+06} & \bftab{9.43290e+06} & 9.47386e+06 & 1.05951e+07 \\
10 &        1 & 7.72013e+06 & 7.72013e+06 &   \bftab{8.37902e+06} & \bftab{8.37902e+06} & 8.54589e+06 & 9.83225e+06 \\
\bottomrule
\end{tabular}
\end{center}
\caption{Heuristic performance for selected datasets. Best upper bounds found are typeset in boldface.}\label{tab:heuristic}
\end{table}

\section{Conclusions}\label{sec:conclusion}
We developed an exact solution algorithm for the minimum sum-of-squares clustering problem (MSSC) using tools from semidefinite programming. 
We use a semidefinite relaxation that exploits three types of valid inequalities in a cutting plane fashion to generate tight lower bounds for the MSSC.

Besides these lower bounds, the semidefinite relaxation also provides a primal solution that can be used for generating data to initialize constrained $k$-means, which is known to be sensitive concerning the starting point. Numerical experiments undoubtedly demonstrate the advantage of using this initialization procedure.

We implemented a branch-and-bound algorithm using the ingredients described above. Our way of branching allows us to decrease the size of the problem while going down the branch-and-bound tree. Notably, the shrinking procedure preserves the structure of the problem which is beneficial for our routine computing the bounds in each node of the branch-and-bound tree. 
Our code is parallelized in two ways: the nodes in the branch-and-bound tree are evaluated in parallel and the bound computation within a node is executed in a multi-threaded MATLAB environment.

The numerical results impressively exhibit the efficiency of our algorithm: we can solve real-world instances up to 4000~data points. To the best of our knowledge, no other exact solution methods can handle generic instances of that size.
Moreover, the dimension of the data points does not influence the performance of our algorithm, we solve instances with more than 20\;000 features.

Our algorithm can be extended to deal with certain constrained versions of sum-of-squares clustering like those with diameter constraints, split constraints, density constraints, or capacity constraints \citep{davidson2005clustering, duong:2017}. This is left for future work. Also, kernel-based clustering is a promising extension that we plan to consider \citep{dhillon2004kernel}.
Finally, we have ideas in mind on how to use our algorithm in a heuristic fashion for obtaining high quality solutions for huge graphs.

\ifJOC
\else

\section*{Acknowledgements}
Parts of this project were carried out during a research stay of the third author at the University Tor Vergata, funded by the University of Rome Tor Vergata Visiting Professor grant 2018.  
Furthermore, this project has received funding from the European Union's Horizon 2020 research and innovation programme under the Marie Sk{\l}odowska-Curie grant 
agreement MINOA No 764759.

We thank Kim-Chuan Toh for bringing our attention to~\cite{JaChayKeil2007} and 
for providing an implementation of the method therein.

\bibliographystyle{plainnat}
\bibliography{bibfile}

{\small
  \vspace*{1ex}\noindent
  Veronica Piccialli,
  \href{mailto:veronica.piccialli@uniroma2.it}{\url{veronica.piccalli@uniroma2.it}},
  University of Rome Tor Vergata, 
  via del Politechnico, 00133 Rome, Italy
  ORCiD: 0000-0002-3357-9608 

  \vspace*{1ex}\noindent
  Antonio M. Sudoso,
  \href{mailto:antonio.maria.sudoso@uniroma2.it}{\url{antonio.maria.sudoso@uniroma2.it}},
  University of Rome Tor Vergata, 
  via del Politecnico, 1, 00133 Rome, Italy
  ORCiD: 0000-0002-2936-9931
  
  \vspace*{1ex}\noindent
  Angelika Wiegele,
  \href{mailto:angelika.wiegele@aau.at}{\url{angelika.wiegele@aau.at}},
  Alpen-Adria-Universität Klagenfurt,
  Universitätsstraße 65--67, 9020 Klagenfurt, Austria, 
  ORCiD:  0000-0003-1670-7951
}

\fi

\end{document}